\documentclass[12pt, a4paper]{amsart}

\usepackage[utf8]{inputenc}
\usepackage{amsthm}
\usepackage{dsfont}
\usepackage{amsmath}
\usepackage{amssymb}
\usepackage{graphicx}
\usepackage{stmaryrd}
\usepackage{commath}
\usepackage{tikz}
\usepackage{tikz-cd}
\usepackage[colorlinks=true,linktocpage=true,pagebackref=false, citecolor=red,linkcolor=blue]{hyperref}
\usepackage[shortlabels]{enumitem}
\usepackage{caption}
\setlist{topsep=-1pt}

\usepackage{geometry}
\newtheorem{tma}{Theorem}[section]
\theoremstyle{definition}
\newtheorem{defi}[tma]{Definition}
\newtheorem{obs}[tma]{Remark}
\newtheorem{ejem}[tma]{Example}
\theoremstyle{plain}
\newtheorem{prop}[tma]{Proposition}
\newtheorem{prodef}[tma]{Proposition and Definition}
\newtheorem{lema}[tma]{Lemma}
\newtheorem{coro}[tma]{Corollary}

\newcommand{\Bd}{\textup{\bf{}Bd}}
\newcommand\Mod[1]{{#1}\textrm{-}\mathbf{Mod}}
\newcommand\Var[1]{\mathbf{Var}_{#1}}
\newcommand\K[1]{\mathrm{K}#1}
\newcommand{\bk}{{\mathbb{K}}}
\newcommand\Varrel[1]{\mathbf{Var}/{#1}}
\newcommand\Unit[1]{\mathds{1}_{#1}}
\newcommand{\lL}{\mathbb{L}}
\newcommand\Span[1]{\text{Span}(#1)}
\newcommand\Aff[2]{\mathrm{AGL}_{#1}(#2)}

\newcommand{\Z}{\mathbb{Z}}
\newcommand{\sS}{\mathbb{S}}
\newcommand{\C}{\mathbb{C}}

\newcommand{\Hom}{\operatorname{Hom}}
\newcommand{\id}{\operatorname{id}}
\newcommand{\Stab}{\operatorname{Stab}}
\newcommand{\SL}{\mathrm{SL}}
\newcommand{\PGL}{\mathrm{PGL}}
\newcommand{\AGL}{\mathrm{AGL}}
\newcommand{\comvar}{\mathcal{C}_2(G)}

\newcommand{\torusinsphere}{%
\begin{tikzpicture}

\colorlet{myaxes}{black!85}
\colorlet{mytorus}{black!80}
\colorlet{myblue}{blue!80!black}
\colorlet{myred}{red!80!black}

\draw[->, very thick, myaxes] (0,0) -- (0, 2.8); 
\node[myaxes] at (0.45, 2.7) {\Large $\beta$};
\draw[->, thin, myaxes] (0,0) -- (3.8, 0); 
\draw[->, thin, myaxes] (0,0) -- (-2.0, -2.8); 
\filldraw[myaxes] (0,0) circle (1pt); 

\draw[thick, mytorus] (0,0) ellipse (3.4 and 2.4);

\draw[thick, mytorus] (-1.5, 0.1) to[bend right=45] (1.5, 0.1);
\draw[thick, mytorus] (-1.2, -0.15) to[bend left=45] (1.2, -0.15);

\draw[thick, dashed, mytorus] (0,0) ellipse (2.6 and 1.6);
\node[mytorus] at (2.2, -1.2) {\Large $\alpha$};

\draw[thick, myblue] (0,0) ellipse (3.2 and 2.2);
\node[myblue] at (2.7, -1.7) {\Large $a$};

\def\rx{1.3}
\def\ry{-1.38}
\draw[thick, myred] (\rx, \ry+1.25) arc (90:-90:0.35 and 1.05);
\node[myred] at (1.7, -0.3) {\Large $b$};
\draw[thick, myred, dashed] (\rx, \ry-0.85) arc (270:90:0.35 and 1.05);

\end{tikzpicture}
}

\newcommand{\vankampen}{%
\begin{tikzpicture}

  \tikzset{
    tube line/.style={very thick, black!85},
    tube dashed/.style={thick, black!85, dashed},
    tube loosely dashed/.style={thick, black!85, loosely dashed},
    red solid/.style={thick, black!80!black, fill=black!70},
    red dashed/.style={thick, black, dashed},
    red shade/.style={fill=black!15},
    grey shade/.style={fill=black!5}
  }

  \begin{scope}
    \draw[tube line] (0,0) ellipse (0.9 and 1.6);
  \end{scope}

  \path[grey shade] (2.8, -1.6) .. controls (3.6, -1.5) and (3.8, 1.2) .. (4.8, 2.48) -- (5.5, 3) .. controls (4.5, 0.8) and (4.1, -1.9) .. (3.5, -2) -- cycle;
  \draw[red dashed] (2.8, -1.6) .. controls (3.6, -1.6) and (3.8, 1.2) .. (4.8, 2.48);
  \draw[red dashed] (5.5, 3) .. controls (4.5, 0.8) and (4.1, -1.9) .. (3.5, -2);

  \draw[tube loosely dashed] (0.0, 1.6) .. controls (2.5, 1.4) and (4, 1.6) .. (5.5, 3);
  \draw[tube line] (0.0, 1.6) .. controls (0.3, 1.6) and (0.6, 1.55) .. (1, 1.55);
  \draw[tube loosely dashed] (0, -1.6) .. controls (2.5, -1.8) and (4, -1.6) .. (5.5, -0.2);
  \draw[tube line] (0.0, -1.6) .. controls (0.3, -1.62) and (0.6, -1.63) .. (0.9, -1.64);
  
  \path[red shade] (2.8, -1.6) .. controls (3.2, -1.45) and (3.6, -1.45) .. (4, -1.23) -- (4, -1.23) .. controls (3.8, -1.78) and (3.65, -1.89) .. (3.5, -2) -- cycle;
  \path[red shade] (0.6, 1.2) .. controls (1.2, 0.9) and (1.8, 0) .. (2.8,-1.6) -- (3.5, -2) .. controls (2.5, -0.1) and (1.9, 1.2) .. (0.8, 1.8) -- cycle;
  
  \draw[red dashed] (3.5, -2) .. controls (2.5, -0.1) and (1.9, 1.2) .. (0.8, 1.8);
  \draw[red dashed] (0.6, 1.2) .. controls (1.2, 0.9) and (1.8, 0) .. (2.8,-1.6);

   \path[red shade] (4.8, 2.48) .. controls (5.2, 3) and (5.3, 3.2) .. (5.5, 3.3) -- (5.5, 2.6) .. controls (5.3, 2.65) and (5.1,2.6 ) .. (4.8, 2.48) -- cycle;
   \node[tube line] at (4.1, 0) {\Large $R$};
  \draw[red dashed] (4.8, 2.48) .. controls (5.2, 3) and (5.3, 3.2) .. (5.5, 3.3);
  \draw[red dashed] (5.5, 2.6) .. controls (5.3, 2.65) and (5.1,2.6 ) .. (4.8, 2.48);

  \draw[tube line] (0, 3) .. controls (2.5, 2.8) and (4.0, 3) .. (5.5, 4.6);
  \draw[tube line] (0, -3) .. controls (2.5, -3.2) and (4.0, -3) .. (5.5, -1.6);

  \draw[tube dashed] (0, 2.2) .. controls (2.5, 2) and (4, 2.2) .. (5.5, 3.8);
  \draw[tube line] (0, 2.19) .. controls (0.4, 2.17) and (0.8, 2.16) .. (1.4, 2.15);
  \draw[tube dashed] (0, -2.2) .. controls (2.5, -2.4) and (4, -2.2) .. (5.5, -0.8);
  \draw[tube line] (0, -2.21) .. controls (0.4, -2.23) and (0.8, -2.24) .. (1.3, -2.25);

  \begin{scope}
    \draw[tube line] (0,0) ellipse (2 and 3);
    \node[tube line] at (-1.8, 2) {\Large $\mathbb{T}$};
    \draw[tube line] (0,0) ellipse (1.4 and 2.2);
    \node[tube line] at (-0.4, 0.8) {\Large $\mathbb{T}^{\dagger}$};
  \end{scope}

\end{tikzpicture}
}

\definecolor{thicktorusfill}{RGB}{248,246,244}
\definecolor{thicktorusline}{RGB}{25,25,25}
\definecolor{knotboundarycolor}{RGB}{138,20,35}

\newcommand{\thickenedtoruscut}{%
\begin{tikzpicture}[
  baseline=(current bounding box.center),
  line cap=round,
  line join=round,
  torusboundary/.style={
    draw=thicktorusline,
    line width=1.05pt
  },
  knotboundary/.style={
    draw=knotboundarycolor,
    line width=1.15pt
  },
  knotleader/.style={
    draw=knotboundarycolor,
    line width=0.55pt
  },
  toruslabel/.style={
    inner sep=1pt,
    text=thicktorusline,
    font=\normalsize
  },
  regionlabel/.style={
    inner sep=1pt,
    text=thicktorusline!60,
    font=\small
  },
  knotlabel/.style={
    inner sep=1pt,
    text=knotboundarycolor,
    font=\footnotesize
  }
]

  \path[fill=thicktorusfill, even odd rule]
    (0,0) circle[radius=2cm]
    (0,0) circle[radius=0.68cm];

  \draw[torusboundary] (0,0) circle[radius=2cm];
  \draw[torusboundary] (0,0) circle[radius=0.68cm];

  \coordinate (freeboundary) at (1.15,0.40);
  \fill[white] (freeboundary) circle[radius=0.30cm];
  \draw[knotboundary] (freeboundary) circle[radius=0.30cm];

  \node[toruslabel] at (-0.91,2.1) {$\mathbb{T}$};
  \node[toruslabel] at (0.2,-0.9) {$\mathbb{T}$};
  \node[regionlabel] at (-1.36,-0.12) {$\Pi_{n,m}$};

  \node[knotlabel, anchor=south east] (freeboundarylabel)
    at (1.5,0.8) {$\partial\nu(K_{n,m})\cong\mathbb{T}$};
  \draw[knotleader]
    (freeboundarylabel.south east)
    to[out=-90,in=55] (1.32,0.63);

  \pgfresetboundingbox
  \path[use as bounding box]
    (-2.34,-2.34) rectangle (2.34,2.34);

\end{tikzpicture}%
}

\title[A TQFT for Representation Varieties of Torus Links]{A Topological Quantum Field Theory for Representation Varieties of Torus Links}


\author{Enrique Aycart-Maldonado}
\address{Departamento de Álgebra, Geometría y Topología,
Facultad de Ciencias Matemáticas,
Universidad Complutense de Madrid,
Plaza de las Ciencias, 3,
28040 Madrid, Spain}
\email{enaycart@ucm.es}

\author{Alejandro Calleja}
\address{Departamento de Álgebra, Geometría y Topología,
Facultad de Ciencias Matemáticas,
Universidad Complutense de Madrid,
Plaza de las Ciencias, 3,
28040 Madrid, Spain}

\address{Instituto de Ciencias Matemáticas, Consejo Superior de Investigaciones Científicas, Nicolás Cabrera, 15, 28049 Madrid, Spain}
\email{alejac03@ucm.es}

\begin{document}

\begin{abstract}
We construct a Topological Quantum Field Theory (TQFT) that computes the virtual classes of $G$-representation varieties associated with the fundamental group of a torus link complement. Following the framework established by González-Prieto, Logares and Muñoz for $G$-representation varieties of surface groups, we introduce a suitable bordism whose iteration enables a recursive topological evaluation of these virtual classes. Via this TQFT, we explicitly compute the virtual class of these varieties for the groups $G=\AGL_1(\C)$ and $G=\SL_2(\C)$. 
\end{abstract}

\maketitle

\section{Introduction}

Let $\Gamma$ be a finitely presented group and let $G$ be an algebraic group over a field $\bk$. It is an easy fact that the set $R_G(\Gamma)=\Hom(\Gamma,G)$ of all representations $\rho:\Gamma\to G$ of $\Gamma$ into $G$ has the structure of an algebraic variety over $\bk$, called the \emph{$G$-representation variety} of $\Gamma$. This variety has a natural action of the group $G$ by conjugation of the representations. Therefore, if $G$ is reductive, we can consider the Geometric Invariant Theory (GIT) quotient $\mathfrak{X}_G(\Gamma)=R_G(\Gamma)\sslash G$, called the $G$-\emph{character variety} of $\Gamma$. 

These varieties have been widely studied in recent years, mainly due to their relation with the non-abelian Hodge correspondence \cite{corlette1988flat,hitchin1987self,simpson1994moduli}. This correspondence relates the moduli space of polystable $G$-Higgs bundles over a smooth complex projective curve $\Sigma$ with the character variety $\mathfrak{X}_G(\pi_1(\Sigma))$ and with the moduli space of reductive flat $G$-connections.

In this setting, motivated by mirror symmetry, Hausel and Thaddeus conjectured \cite{hauselthaddeus2003} an equality between appropriately twisted stringy $E$-polynomials of moduli spaces associated with the Langlands dual groups $\SL_n(\C)$ and $\PGL_n(\C)$, in the coprime-degree setting. This conjecture was proved in arbitrary rank by Groechenig, Wyss, and Ziegler \cite{groechenigwyssziegler2020} for the corresponding moduli spaces of Higgs bundles.

These developments have motivated several techniques to compute $E$-polynomials and other motivic invariants of character and representation varieties: the \emph{arithmetic method}, developed by Hausel and Rodríguez-Villegas \cite{hauselvillegas2008katz}, based on counting points over finite fields; the \emph{geometric method}, introduced in \cite{lomune2013surfaces}, which consists in stratifying the varieties into simpler pieces; and the \emph{quantum method}, introduced in \cite{gonzalezprieto2018TQFTcharvar,GONZALEZPRIETOlaxmonoTQFTrepvars} by González-Prieto, Logares and Muñoz, where the authors constructed a lax monoidal Topological Quantum Field Theory (TQFT for short) computing the Grothendieck motive or virtual class of representation varieties in the Grothendieck ring of algebraic varieties $\K{\Var{\bk}}$. Recall that the Grothendieck motive is universal for additive and multiplicative invariants, in the sense that any such invariant factorizes through $\K{\Var{\bk}}$.

A very interesting application of character varieties appears in knot theory. Recall that the group of a knot (or a link) $K \subseteq \sS^3$ is defined as the fundamental group $\pi_1(\sS^3-K)$ of its complement in the $3$-dimensional sphere. This group is finitely generated, so one can study its character and representation varieties. This has been classically done due to its relation with classical knot invariants, something that can be appreciated in the works of de Rham \cite{deRham1968introduction}, where a relation is shown between the $\SL_2(\C)$-character variety of the knot and its Alexander polynomial; or of Cooper et al.\ \cite{cooper1994Apoly}, where they define one of the most important knot invariants, the $A$-polynomial, in terms of its $\SL_2(\C)$-character variety. 

Nevertheless, in recent years, character varieties of knots have become of interest in their own right, mainly motivated by the theory of character varieties of surfaces. In this sense,  Florentino, Nozad and Zamora proved in \cite{florentino2021serre} a mirror-symmetric statement in the sense of \cite{hauselvillegas2008katz}, by showing the equality of $E$-polynomials of the $\SL_n(\C)$ and $\PGL_n(\C)$-character varieties of trivial links, using combinatorial methods. This behavior has not yet been observed in any other knot or link, mainly due to the difficulty of computing this invariant. In fact, very few examples of $E$-polynomials of these varieties have been computed, since the only knots for which this invariant is known are the figure-eight knot, whose $\SL_3(\C)$-character variety was studied in \cite{heusener2016sl3figureeight}; and the torus knots $K_{n,m}$, which were studied in \cite{deAmorin2025motive,GONZALEZPRIETOtorusknotsSL4,martinezmartinez2012su2torusknots,munoz2009sl2torusknot,munoz2016sl3torusknots}. The computation of all these examples was done by using the aforementioned geometric method. 

The situation with links is considerably less explored, since the only case studied is that of torus links $L_{n,m}^d$. The $E$-polynomials of the corresponding character varieties for $\SL_2(\C)$ and $\SL_3(\C)$ were recently computed in \cite{gonzalez2024toruslinks}. In that paper, the authors provide a stratification of the representation varieties in terms of the stabilizers of the conjugation action, which enables the construction of the character variety. Crucially, this work establishes evidence that the $E$-polynomial of the character varieties for this class of links can be computed via a TQFT, adapting the framework developed for surface groups \cite{gonzalezprieto2018TQFTcharvar,GONZALEZPRIETOlaxmonoTQFTrepvars}. In these foundational works of the quantum method, the authors construct a TQFT through the composition of two functors: the geometrization and quantization functors. The former assigns to a bordism the $G$-representation variety of its fundamental group, while the latter maps varieties to modules over the Grothendieck ring of varieties, and sends spans to module homomorphisms. Following this technique, a recursive topological decomposition is applied to analyze the representation varieties of genus $g$ surface groups, assembling them from elementary pieces of genus 1. 

Following this approach, we will show that the TQFT constructed in \cite{gonzalezprieto2018TQFTcharvar,GONZALEZPRIETOlaxmonoTQFTrepvars} can be used to compute the Grothendieck motive of representation varieties of torus links by giving an appropriate description of the manifold $\sS^3-L_{n,m}^d$ as a composition of bordisms. Moreover, we will show that this TQFT generalizes the stratification given in \cite{gonzalez2024toruslinks}, which suggests that it can also be applied to the character variety. We will use this method to compute the motive of the representation varieties of torus links for the algebraic groups $\AGL_1(\C)$ and $\SL_2(\C)$, obtaining the following result:

\begin{tma}
    Let $R_G(L_{n,m}^d)$ be the $G$-representation variety of the torus link $L_{n,m}^d$ with $d, n, m \geq 1$, $n$ and $m$ coprime. Then, for $G = \AGL_1(\C)$ and $G=\SL_2(\C)$ we have that the virtual classes of the respective representation varieties are
    \begin{align*}
    \left[R_{\mathrm{AGL}_1(\C)}(L_{n,m}^d)\right] &= \lL(\lL-1-nm)(\lL-1)^{d-1} + (\lL-1)\lL^{d-1} \\
    & + (1+(m-1)\lL)(1+(n-1)\lL)(\lL(\lL-1))^{d-1},
\end{align*}
and
   \begin{align*}
[R_{\SL_2(\C)}(L_{n,m}^d)] &= \left(1 + \frac{m-1}{2}(\lL^2+\lL)\right) \left(2 + (n-1)(\lL^2+\lL)\right) (\lL^3-\lL)^{d-1} \\
&\quad + \left(\lL^3 - (mn+1)\lL^2 - mn\lL\right) (\lL+1)^{d-1} \\
&\quad + \frac{1}{2} (\lL - 1 - 2mn)(\lL^2+\lL) \left( (\lL-1)^{d-1} - (\lL+1)^{d-1} \right) \\
&\quad + 2(\lL^2-1)(2\lL)^{d-1}.
\end{align*}
    where $\lL=[\mathbb{A}_{\C}^1]$ is the Lefschetz motive.
\end{tma}

Moreover, we recover some results found in the literature as direct corollaries of the application of the quantum method. In particular, we recover the virtual class of the $\AGL_1(\C)$-representation variety of a torus knot \cite{gplogmunlowrankAGL1} and the $E$-polynomial of the $\SL_2(\C)$-representation variety of a twisted Hopf link \cite{gonzalez2023hopflinks}. Although the explicit computations in this paper focus on
$\AGL_1(\C)$ and $\SL_2(\C)$, the techniques developed here are general and can be applied to other algebraic groups. In particular, the method does not require reductivity and provides a framework for computing virtual classes of torus link representation varieties for arbitrary affine algebraic groups, provided that the corresponding stratification into stabilizer pairs is finite.

This work is structured as follows. In Section \ref{sec:rep_varieties} we introduce the representation variety of a link, focusing on the case of torus links. We also give a description of these links that will later allow us to express the corresponding manifold as a composition of bordisms. In Section \ref{sec:Quantum_Method}, we review the construction of the TQFT introduced in \cite{gonzalezprieto2018TQFTcharvar,GONZALEZPRIETOlaxmonoTQFTrepvars}, adapting it to the setting of \emph{horizontally convergent bordisms with free boundary}, which allows us to simplify the computations. In Section \ref{sec:QuantizationVirtualClasses}, we use the description of a torus link $L_{n,m}^d$ introduced in Section \ref{sec:rep_varieties} to obtain a decomposition of $\sS^3-L_{n,m}^d$ into bordisms and compute the corresponding TQFT operators for any algebraic group. Finally, in Sections \ref{sec:AGL1} and \ref{sec:SL2}, we apply the computations of Section \ref{sec:QuantizationVirtualClasses} to the groups $\AGL_1(\C)$ and $\SL_2(\C)$, respectively, obtaining explicit formulas for the virtual classes of the corresponding representation varieties.

\subsection*{Acknowledgements}
The authors would like to express their sincere gratitude to their PhD supervisor, Ángel González-Prieto, for his invaluable guidance, insightful suggestions, and support throughout the preparation of this work. They also wish to thank Francisco Javier Martínez-Aguinaga, Lucas de Amorin, and Carlos Florentino for their willingness to discuss this work and for their valuable comments and suggestions.

This work was supported by projects PID2021-124440NB-I00 and PID2024-156578NB-I00, funded by MICIU/AEI/10.13039/501100011033 and by the European Regional Development Fund (ERDF, EU). Alejandro Calleja is also partially supported by CEX2023-001347-S, funded by MICIU/AEI/10.13039/501100011033.

\section{Representation Varieties of Torus Links}\label{sec:rep_varieties}

In this section, we introduce the links and representation varieties that will be studied throughout the paper. We first recall the definition of the representation and character varieties associated with a link group. We then specialize to torus links and recall a presentation of their fundamental groups. In order to use the nested solid-torus decomposition underlying our constructions, we introduce nested torus links via a radial compression of $\sS^3$. Finally, we prove that every standard torus link is ambiently isotopic to a nested link consisting of parallel copies of the corresponding torus knot. This alternative description will be particularly useful in the subsequent decomposition of their complements.

Let $L\subseteq \sS^3$ be a link, and consider its group
$$\Gamma_L:=\pi_1(\sS^3-L)=\langle\gamma_1,\dots,\gamma_r \ | \ R_j(\gamma_1,\ldots,\gamma_r), \ j=1,\ldots,s\rangle,$$
which is finitely presented thanks to Wirtinger's algorithm. Throughout the paper, $\bk$ denotes an algebraically closed field of characteristic zero. Let $G$ be an affine algebraic group over $\bk$. We define the \emph{$G$-representation variety} of the link $L$ as the set 
$$R_G(L):=\Hom(\Gamma_L,G)=\{(g_1,\dots,g_r)\in G^r \ | \ R_j(g_1,\dots,g_r)=e, \ j=1,\dots, s\},$$
which is an affine variety over $\bk$ (we do not require varieties to be irreducible).

In the case when $G$ is reductive, one can also define the \emph{$G$-character variety} of the link as the GIT quotient 
$$\mathfrak{X}_G(L):=R_G(L)\sslash G,$$
where $G$ acts on the representation variety by conjugation. Recall that given an affine variety $Z=\mathrm{Spec}(R)$ for a finitely generated $\bk$-algebra $R$, and an action of a reductive group $G$ on $Z$, the Geometric Invariant Theory (GIT for short) quotient of $Z$ is defined as 
$$Z\sslash G:=\mathrm{Spec}(R^G),$$
where $R^G$ is the algebra of $G$-invariant elements of $R$.

\begin{figure}[ht]
    \centering
    \begin{minipage}{0.48\textwidth}        
        \centering
        \includegraphics[width=0.8\linewidth]{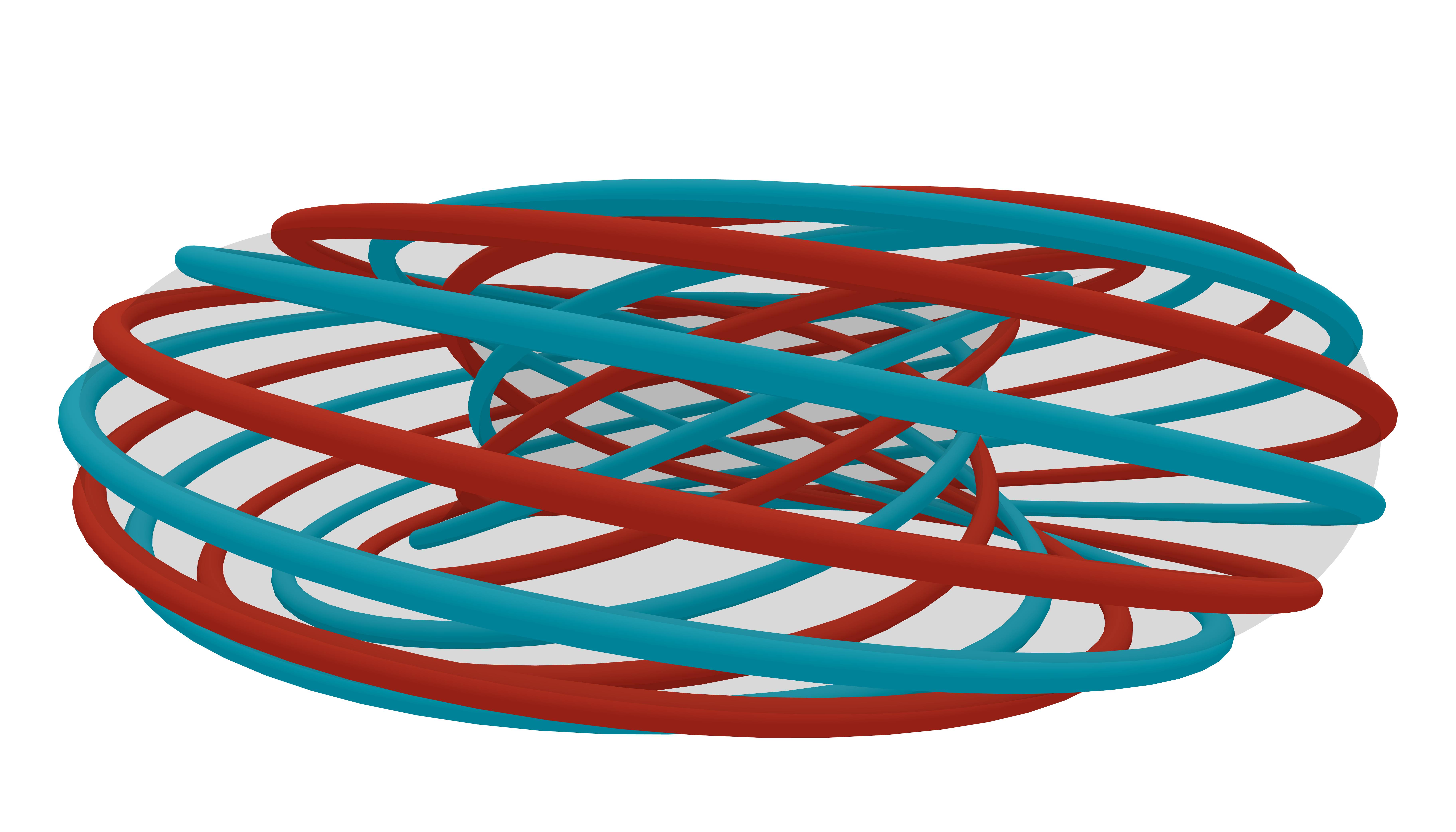}
    \end{minipage}
    \hfill
    \begin{minipage}{0.48\textwidth}
        \centering
        \includegraphics[width=1\linewidth]{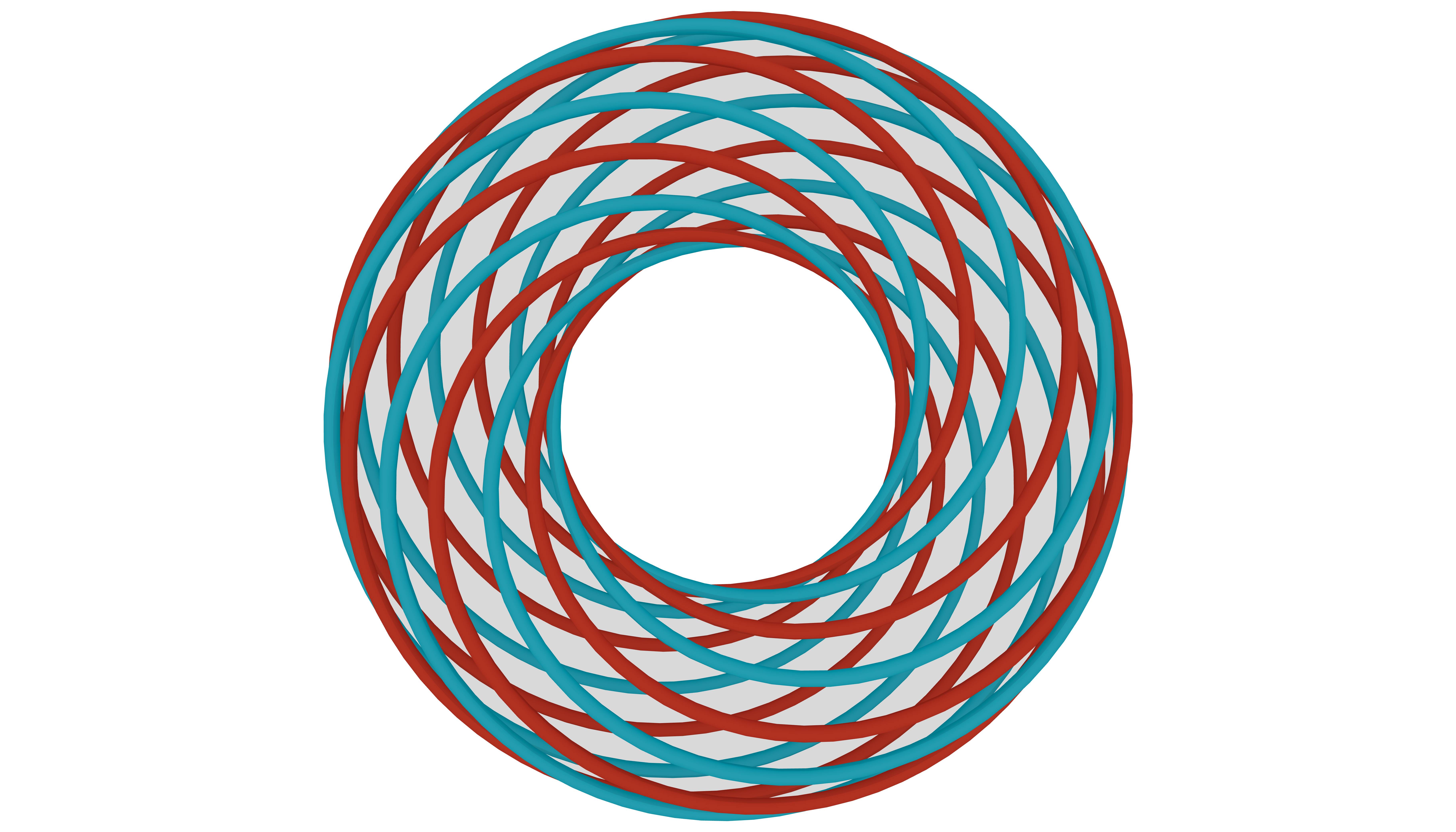}
    \end{minipage}

    \caption{Frontal and top views of the torus link $L^2_{5,7}$, that is, the two-component $(5,7)$-torus link.}
    \label{fig:front-upper}
\end{figure}

In this paper, we will focus on the case of \emph{torus links}. Recall that given two coprime integers $n,m \geq 1$, the torus knot $K_{n,m}$ is the image of the embedding
\begin{equation}\label{embtorusknots}
    \mathbb{S}^1 \longrightarrow \sS^3\subseteq\C^2, \qquad z\longmapsto \tfrac{1}{\sqrt{2}}\,(z^n,z^m),
\end{equation}
or, equivalently, the image of a line with rational slope $n/m$ in the square representation of the torus. More generally, one can define the torus link $L_{n,m}^d$ as the image of $d$ parallel lines in the square representation. This yields a link with $d$ components, which can equivalently be viewed as the $(dn,dm)$-torus link. Two views of the two-component $(5,7)$-torus link $L_{5,7}^2$ are shown in Figure \ref{fig:front-upper}.

It was proven in \cite{argyres2019fundamentalgrouptoruslinks} that the group of this link is
$$\Gamma_{n,m}^d:=\Gamma_{L_{n,m}^d}=\langle a,b,f_1,\ldots,f_{d} \ | \ f_1=e,a^nf_j=f_jb^m, j=1,\dots,d\rangle,$$
where $e$ denotes the identity element in $\Gamma_{n,m}^d$.

For the purposes of this work, we will also need to introduce the notion of \emph{nested torus links} (see also \cite{argyres2019fundamentalgrouptoruslinks}). Let $\sS^3$ be parameterized by the hyperspherical coordinates $(\chi, \theta, \psi) \mapsto \left(\cos( \chi) e^{i\theta},\sin( \chi) e^{i\psi}\right)$, with $0 \leq \chi \leq \pi/2$, which naturally induce a decomposition into two solid tori $\{(\chi, \theta, \psi) \ | \ \chi < \pi/4\}$ and $\{(\chi, \theta, \psi) \ | \ \chi > \pi/4\}$, henceforth denoted as $\mathbf{T}^\dagger$ and $\mathbf{T}$, respectively. Both solid tori are bounded by the auxiliary 2-torus $T = \{\chi = \pi/4\}$. We define a radial compression homeomorphism $\phi: \sS^3 \to \sS^3$ given by $$\phi(\chi, \theta, \psi) = \left(\frac{2}{\pi}\chi^2, \theta, \psi\right),$$ which maps the inner solid torus $\overline{\mathbf{T}^\dagger}$ into its interior. We use this map to inductively place a link inside the auxiliary torus of another link.

\begin{defi}
    A nested torus link of depth $r$, denoted $N_{(n_1,m_1),\ldots,(n_r,m_r)}^{d_1,\ldots,d_r}$, is defined inductively. For the base case $r=1$, it is a standard torus link:$$N_{(n_1,m_1)}^{d_1} := L_{n_1,m_1}^{d_1}.$$For $1 \le k < r$, the depth-$(k+1)$ nested link is the disjoint union of a standard torus link on $T$ and the image of the depth-$k$ nested link under the compression map $\phi$:$$N_{(n_1,m_1),\ldots,(n_{k+1},m_{k+1})}^{d_1,\ldots,d_{k+1}} := L_{n_{k+1},m_{k+1}}^{d_{k+1}} \cup \phi\left(N_{(n_1,m_1),\ldots,(n_k,m_k)}^{d_1,\ldots,d_k}\right).$$
    Equivalently, this gives the explicit finite union:$$N_{(n_1,m_1),\ldots,(n_r,m_r)}^{d_1,\ldots,d_r} = \bigcup_{i=1}^{r} \phi^{r-i}\left(L_{n_i,m_i}^{d_i}\right).$$
\end{defi}

\begin{figure}[ht]
    \centering
    \includegraphics[width=0.7\linewidth]{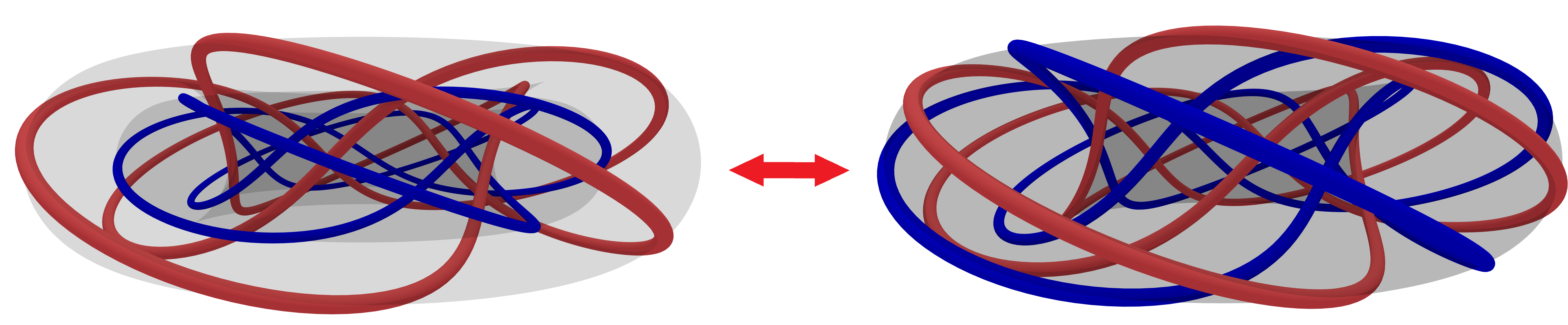}
    \caption{Schematic illustration of the isotopy described in Proposition \ref{prop:toruslinksarenested} for an $(n,m)$-nested torus link of depth $d=2$.}
    \label{fig:isotopy}
\end{figure}

An important observation is that standard torus links admit an equivalent description in terms of nested torus knots. More precisely, every torus link $L_{n,m}^d$ is ambiently isotopic to a nested link of depth $d$, consisting of $d$ copies of the $(n,m)$-torus knot. A schematic illustration of the aforementioned isotopy is shown in Figure \ref{fig:isotopy}.

\begin{prop}\label{prop:toruslinksarenested}
    There is an ambient isotopy between the torus link $L_{n,m}^d$ and the link of $(n,m)$-nested torus knots of depth $d$.
\end{prop}

\begin{proof}
The $d$ components of $L_{n,m}^d$ are parallel copies of the $(n,m)$-torus knot lying on the auxiliary torus $T=\{\chi=\pi/4\}$. We fix one of the components. The remaining ones may be isotoped successively. We consider the compression map $\phi$ and move the $j$-th component onto the torus $\phi^{j-1}(T)$ with $j=1,...,d$. It is clear that $\phi$ preserves the slope and, since these tori are pairwise disjoint and the compression does not modify the coordinates $(\theta,\psi)$, in which the components are already disjoint, the components remain disjoint throughout the compression. After the process the result is $$\bigcup_{j=0}^{d-1}\phi^j(L^1_{n,m})=N^{\scriptstyle\overbrace{1,\ldots,1}^{d\text{ times}}}_{\scriptstyle\underbrace{(n,m),\ldots,(n,m)}_{d\text{ times}}},$$
which is a depth-$d$ nested link of $(n,m)$-torus knots. Thus we obtain an isotopy of the link and it extends to an ambient isotopy of $\sS^3$. 
\end{proof}

\section{Quantum Method}\label{sec:Quantum_Method} 

This section introduces the construction of the Topological Quantum Field Theory (TQFT) that computes the virtual class of a $G$-representation variety. We shall consider the framework developed by González-Prieto, Logares and Muñoz \cite{gonzalezprieto2018TQFTcharvar, GONZALEZPRIETOlaxmonoTQFTrepvars}, and adapt it to bordisms with free boundary.

\subsection{The category of bordisms with free boundary}

For the purposes of this work, we introduce a larger type of bordisms that provides the appropriate framework for the constructions considered in this paper.

\begin{defi}[Bordisms with free boundary]
Let $M_-$ and $M_+$ be closed $(n-1)$-manifolds, possibly empty. An \emph{$n$-dimensional bordism with free boundary} from $M_-$ to $M_+$ is a compact smooth $n$-manifold $W$, possibly disconnected, together with a decomposition of its boundary
$$\partial W=\partial^-_{\mathrm h}W\sqcup\partial^+_{\mathrm h}W\sqcup\partial_{\mathrm f}W$$
into disjoint unions of connected components, and diffeomorphisms
$$\iota_-:M_-\xrightarrow{\sim}\partial^-_{\mathrm h}W, \qquad \iota_+:M_+\xrightarrow{\sim}\partial^+_{\mathrm h}W.$$

The submanifolds $\partial^-_{\mathrm h}W$ and $\partial^+_{\mathrm h}W$ are called the \emph{incoming} and \emph{outgoing horizontal boundaries}, respectively, whereas $\partial_{\mathrm f}W$ is called the \emph{free boundary}. We call the diffeomorphisms $\iota_-$ and $\iota_+$ the \emph{boundary parametrizations}.

Two bordisms $(W,\iota_-,\iota_+),(W',\iota'_-,\iota'_+)$ from $M_-$ to $M_+$ are declared to be \emph{equivalent} if there exists a diffeomorphism $F:W\xrightarrow{\sim}W'$ such that $F(\partial_\mathrm fW)=\partial_\mathrm f W'$ and the following diagram commutes:
    \begin{equation}\label{diag:diagram_def_bdp}
\begin{tikzcd}[column sep=large, row sep=large]
&
W \arrow[dd, "F", "\sim"']
&
\\
M_-
  \arrow[ur, "\iota_-"]
  \arrow[dr, "\iota'_- "']
&
&
M_+
  \arrow[ul, "\iota_+"']
  \arrow[dl, "\iota'_+"]
\\
&
W' &
\end{tikzcd}
\end{equation}
\end{defi}

\begin{ejem}
Consider the torus $\mathbb{T}=\sS^1\times \sS^1$ and the cylinder $\mathbb{T}\times[0,1]$. Let $B^3$ be a closed $3$-ball smoothly embedded in the interior $\mathbb{T}\times(0,1)$, and define
$$
W=(\mathbb{T}\times[0,1])\setminus\operatorname{int}(B^3).
$$
Geometrically, $W$ is obtained from a cylinder over the torus by removing an open ball from its interior. Its boundary consists of three connected components: the two tori at the ends of the cylinder and the sphere surrounding the removed ball.

Thus, $W$ is a $3$-dimensional bordism with free boundary from $\mathbb{T}$ to $\mathbb{T}$, with
$$
\partial^-_{\mathrm h}W=\mathbb{T}\times\{0\},\qquad
\partial^+_{\mathrm h}W=\mathbb{T}\times\{1\},\qquad
\partial_{\mathrm f}W=\partial B^3\cong \sS^2.
$$
The boundary parametrizations are given by $\iota_-(x)=(x,0)$ and $\iota_+(x)=(x,1)$. Notice that the spherical boundary component is not parametrized, since it belongs to the free boundary.
\end{ejem}

\begin{defi}
The category $\Bd^{\mathrm{free}}_n$ has closed $(n-1)$-manifolds, including $\varnothing$, as objects. Its morphisms $M_-\to M_+$ are equivalence classes of $n$-dimensional bordisms with free boundary as above.

If $W:M_0\to M_1$ and $V:M_1\to M_2$ are morphisms, their composite is obtained by gluing $\partial^+_{\mathrm h}W$ to $\partial^-_{\mathrm h}V$ through the boundary parametrizations. The free-boundary components are not glued.
\end{defi}

In what follows, we restrict our attention to a particular class of bordisms with free boundary.

\begin{defi}
A morphism $W:M_-\to M_+$ in $\Bd^{\mathrm{free}}_n$ is \emph{horizontally convergent} if every connected component $C$ of $W$ satisfies
$$\#\pi_0(C\cap\partial^-_{\mathrm h}W)\leq 1, \qquad \#\pi_0(C\cap\partial^+_{\mathrm h}W)\leq 1.$$
There is no restriction on $C\cap\partial_{\mathrm f}W$.
\end{defi}

\begin{defi}
We denote by $\overline{\Bd}_n^\mathrm{free}$ the \emph{category of horizontally convergent bordisms with free boundary}, defined as follows:
\begin{itemize}
\item \emph{Objects}: The objects are closed oriented $(n-1)$-manifolds equipped with one marked point in each connected component, together with the empty manifold $\varnothing$.

\item \emph{Morphisms}: A morphism is a horizontally convergent bordism together with the following additional data. For each connected component having both nonempty incoming and outgoing horizontal boundary, we specify a homotopy class, relative to the endpoints, of paths joining the corresponding marked points. Connected components with at most one nonempty horizontal boundary carry no additional path data.
\end{itemize}

The equivalence relation on bordisms is required to preserve the prescribed homotopy classes of paths. Composition is given by gluing bordisms along their common horizontal boundary and, whenever the resulting connected component has both nonempty incoming and outgoing horizontal boundary, by concatenating the corresponding path classes. No path data are retained on components with at most one nonempty horizontal boundary. The identity morphism is the identity cylinder, equipped on each connected component with the vertical path through the corresponding marked point. We shall suppress the marked points and path data from the notation whenever no confusion can arise.
Under disjoint union, the category of horizontally convergent bordisms with free boundary, $( \overline{\Bd}_n^\mathrm{free}, \sqcup, \varnothing )$ is a symmetric monoidal category.
\end{defi}

The role of the marked points and chosen path classes will become clear in the definition of the geometrization functor.

Each connected component of a horizontally convergent bordism has one of the four \emph{horizontal} types
$$
M\longrightarrow P,\qquad
M\longrightarrow\varnothing,\qquad
\varnothing\longrightarrow P,\qquad
\varnothing\longrightarrow\varnothing,
$$
where $M$ and $P$ are nonempty connected closed $(n-1)$-manifolds, representing the incoming and outgoing horizontal boundaries, respectively. Each component may additionally have free-boundary components. For bordisms of the type $W:\varnothing\to\varnothing$, the term \emph{bordism with empty horizontal boundary} is preferable here to \emph{closed manifold}, because the free boundary $\partial_{\mathrm f}W$ may be nonempty.
The importance of this restriction will be explained in Remark \ref{obs:no-basepoints}.

\subsection{The geometrization functor}

For the purposes of the construction, we first recall the definition of the category of spans, which provides the natural target for the geometrization functor.

\begin{defi}
Let $\mathcal C$ be a category admitting fiber products and finite products. The category of spans $\Span{\mathcal C}$ is defined as follows:
\begin{itemize}
    \item \emph{Objects:} The objects of $\mathcal C$.

    \item \emph{Morphisms:} A morphism from $X$ to $Y$ is an isomorphism class of diagrams in $\mathcal C$ of the form
    $$
    X \xleftarrow{f} Z \xrightarrow{g} Y,
    $$
    called \emph{spans}. Two spans are isomorphic if there exists an isomorphism between their middle objects commuting with both arrows.
\end{itemize}
The composition of two spans
$$
X \xleftarrow{f} Z \xrightarrow{g} Y,
\qquad
Y \xleftarrow{h} Z' \xrightarrow{k} T
$$
is defined by the fiber product
$$
X \longleftarrow Z\times_Y Z' \longrightarrow T.
$$
The identity morphism of $X$ is the span $X\xleftarrow{\mathrm{id}}X\xrightarrow{\mathrm{id}}X$.

Moreover, $\Span{\mathcal C}$ is a symmetric monoidal category with monoidal product induced by the Cartesian product in $\mathcal C$ and unit given by its terminal object.
\end{defi}

Taking $\mathcal C=\Var{\bk}$, the category of algebraic varieties over $\bk$, we now construct the geometrization functor.

Let $X$ be a connected compact manifold, and choose a basepoint $x\in X$. We write $\pi_1(X)$ for $\pi_1(X,x)$, suppressing the basepoint from the notation. Let $G$ be an affine algebraic group. The $G$-representation variety of $X$ is defined by
$$R_G(X):=\Hom\bigl(\pi_1(X),G\bigr).$$
Since $\pi_1(X)$ is finitely presented, $R_G(X)$ carries a natural structure of an affine algebraic variety. 

For a possibly disconnected compact manifold $X$, we define 
$$R_G(X):=\prod_{C\in\pi_0(X)} R_G(C),$$
where the marked point is used whenever one is specified, and otherwise a basepoint is chosen in each connected component $C$.

\begin{tma}\label{tma:geometrization}
The assignment
$$
\mathcal F_G:\overline{\Bd}_n^\mathrm{free}\longrightarrow \operatorname{Span}({\Var{\bk}}),
$$
defined on objects by $\mathcal F_G(M)=R_G(M)$ and on bordisms $W:M\to N$ by
$$
\mathcal F_G(W)=
\left(
R_G(M)\xleftarrow{\ \iota_-^*\ }R_G(W)
\xrightarrow{\ \iota_+^*\ }R_G(N)
\right),
$$
where the morphisms are induced by the horizontal boundary parametrizations, is a well-defined strong monoidal functor.
\end{tma}

\begin{proof}
Let $W : M\to N$ be a bordism in $\overline{\Bd}_n^\mathrm{free}$. In a connected component with both horizontal boundaries, we use the incoming marked point as basepoint and the specified path to transport the outgoing boundary loops to this point. In a component with just one horizontal boundary, we choose any basepoint. The boundary inclusions then induce homomorphisms of fundamental groups and hence restriction morphisms. Thus,
$$\mathcal{F}_G(W)= \left( R_G(M)\xleftarrow{\ \iota_{\mathrm{-}}^*\ }R_G(W) \xrightarrow{\ \iota_{\mathrm{+}}^*\ }R_G(N) \right).$$

These spans are well defined up to isomorphism. Indeed, an equivalence of bordisms preserves the marked points and paths, while changing the basepoint (and the path to the marked point) in a component with at most one horizontal boundary only changes its apex by an isomorphism.

To verify compatibility with composition, let $W:M\to N$ and $V:N\to P$. Each pair of connected components that is glued meets along exactly one connected component of $N$. Seifert--van Kampen's theorem at its marked point, with basepoint transport along the specified paths, gives an isomorphism
$$R_G(W\cup_N V)\cong R_G(W)\times_{R_G(N)}R_G(V),$$
compatible with the restriction maps to $R_G(M)$ and $R_G(P)$. For a component with no remaining horizontal boundary only an isomorphism of the apex is required. Taking products over connected components proves compatibility with composition.

The distinguished paths in the identity cylinders are the vertical paths through the marked points, hence the corresponding spans are isomorphic to the identity spans.

Finally, since
$$
R_G(M\sqcup N)\cong R_G(M)\times R_G(N)
$$ and $R_G(\varnothing)=\star$, the geometrization functor $\mathcal F_G$ is strong monoidal with respect to disjoint union and the Cartesian product of varieties.
\end{proof}

Observe that the free boundary is not a third leg of the span, as its loops and their relations are encoded in the representation variety $R_G(W)$, hence in the apex of the span.

\subsection{The Grothendieck ring of varieties}
\label{sec:groth_ring}

In order to properly define the virtual class of an algebraic variety, we shall introduce the construction of the Grothendieck ring of algebraic varieties.

Recall that $\bk$ is an algebraically closed field of characteristic zero. By an \emph{algebraic variety} over $\bk$, we mean a reduced quasi-projective scheme of finite type over $\bk$.

Let $S$ be an algebraic variety over $\bk$. We denote by $\Varrel{S}$ the category of algebraic varieties over $S$. The objects of this category are regular morphisms $X\to S$,
called \emph{$S$-varieties}, and morphisms are given by
regular maps $X\to Y$ preserving the base projections,
called \emph{$S$-morphisms}. In particular, taking $S=\operatorname{Spec}(\bk)$, we recover the category $\Var{\bk}$ of algebraic varieties over $\bk$.

\begin{defi}
The \emph{relative Grothendieck ring of varieties over $S$}, $\K{\Varrel{S}}$, is the commutative ring generated by isomorphism classes $[X\to S]_S$ of objects in $\Varrel{S}$. For a morphism $h:X\to S$, we also denote its class by $[h:X\to S]_S$, or simply $[h:X\to S]$ or $[h]$ whenever the base is clear from the context.

These classes are subject to the scissor relation: for every closed subvariety $Z\subseteq X$,
$$
[X\to S]_S=[Z\to S]_S+[X\setminus Z\to S]_S.
$$
Multiplication is defined by fiber products over $S$:
$$
[X\to S]_S\cdot[Y\to S]_S
=[X\times_S Y\to S]_S,
$$
with unit $[S\xrightarrow{\mathrm{id}}S]_S$. Here, fiber products are understood with their reduced scheme structure.

In particular, taking $S=\operatorname{Spec}(\bk)$, we recover the \emph{Grothendieck ring of varieties over $\bk$}, denoted by $\K{\Var{\bk}}$.
\end{defi}

The class of the affine line $\mathbb{L} = [\mathbb{A}^1_\bk]$ is called the \emph{Lefschetz motive}.

A natural question in this setting is whether two non-isomorphic algebraic varieties can have the same virtual class. The following known result provides a sufficient condition.

\begin{prop} \label{prop:equal_virtual_class}
   Let $f:X \to Y$ be a regular bijection of algebraic varieties. Then $[X]=[Y]$ in $\K\Var{\bk}$. More generally, if $X$ and $Y$ are varieties over a variety $S$, and $f$ preserves the base projections, then $[X\to S]=[Y\to S]$ in $\K\Varrel{S}$.
\end{prop}

\begin{proof}
The argument follows \cite[Lemma 2.1 and Corollary 2.4]{Beke2017}, with the relative statement included here for completeness.

First assume $X$ and $Y$ are irreducible. The reducible case follows by applying the same argument to finitely many irreducible strata.

Since $f$ is a regular bijection of algebraic varieties, it is surjective and quasi-finite (every fiber is finite). Thus, it induces a finite extension of function fields $f^*:\bk(Y)\hookrightarrow\bk(X)$, which is separable since $\operatorname{char}(\bk)=0$. The degree of this extension equals the number of points in a generic geometric fiber. Since $f$ is bijective, this degree is one, and hence $f$ is birational.

Apply now Zariski's Main Theorem \cite[Lemma 37.43.2, Tag 02LR]{stacks-project}. Over the normal locus of $Y$, the quasi-finite birational morphism $f$ is an open immersion, and since it is surjective, it is an isomorphism. In particular, there exists a nonempty open subset $U\subseteq Y$ such that $f^{-1}(U)\cong U$.

Consider $Z=Y\setminus U$. The restriction $f^{-1}(Z)\to Z$ is bijective, while $\dim(Z)<\dim(Y)$. By induction on the dimension, it follows that
$$
[X]=[f^{-1}(U)]+[f^{-1}(Z)]
   =[U]+[Z]=[Y],
$$
concluding the result.

Finally, if $f$ preserves the projections to $S$, the same argument applies, since all the resulting isomorphisms also preserve the projections to $S$. Therefore, $[X\to S]=[Y\to S]$.
\end{proof}

In Section \ref{sec:SL2}, we illustrate the quantum method developed in this paper through explicit computations. To this end, we briefly recall the equivariant virtual class formalism, following \cite[Section~3.6]{vogelthesis}.

Let $\mu_r$ denote the group of $r$-th roots of unity and let $R_{\mathbb Q}(\mu_r)$ be its rational representation ring. The following result follows from \cite[Definition~3.6.12 and Example~3.6.14]{vogelthesis}:

\begin{prodef}\label{prop:equivariant-virtual-class}
Let $X$ be an algebraic variety endowed with an algebraic action of $\mu_r$. There exists a unique class
$$[X]^{\mu_r}\in \K{\Var{\C}}\otimes R_{\mathbb Q}(\mu_r)\otimes\mathbb Q$$ such that, for every subgroup $H\leq \mu_r$,
$$\left\langle 1_H,\operatorname{Res}^{\mu_r}_{H}[X]^{\mu_r}\right\rangle=[X/H],$$
where $1_H$ denotes the trivial representation of $H$, $\operatorname{Res}^{\mu_r}_{H}:R_{\mathbb Q}(\mu_r) \to R_{\mathbb Q}(H)$ is the restriction of representations, and $\langle-,-\rangle$ is the standard inner product of characters, extended linearly to the Grothendieck ring of varieties. 

The element $[X]^{\mu_r}$ is called the \emph{$\mu_r$-equivariant virtual class} of $X$.
\end{prodef}

The multiplicativity of the equivariant virtual class is nontrivial, and it is not satisfied in general. Indeed, equivariant virtual classes need not satisfy $$[X\times Y]^{\mu_r}=[X]^{\mu_r}[Y]^{\mu_r},$$ where $\mu_r$ acts diagonally on $X\times Y$, even for $\mu_r=\mu_2$ (see \cite[Example~3.6.17]{vogelthesis}). 

Recall that $\K\Var{\C}^{\mu_r}$ denotes the Grothendieck ring of algebraic varieties endowed with an algebraic $\mu_r$-action, generated by their equivariant isomorphism classes modulo the relations $[X]=[Z]+[X\setminus Z]$ for $\mu_r$-invariant closed subvarieties $Z\subseteq X$. Multiplication is induced by the Cartesian product with the diagonal action.

Following \cite[Lemma~3.6.18]{vogelthesis}, we consider the multiplicative subalgebra
$$V_{\mu_r}:=\left\{[X]\in \K\Var{\C}^{\mu_r}\ \middle|\ [X\times Y]^{\mu_r}=[X]^{\mu_r}[Y]^{\mu_r} \text{ for every }\mu_r\text{-variety }Y\right\},$$
Thus, in the applications below, the multiplicativity hypothesis means that one of the factors is known to belong to $V_{\mu_r}$. In that case the product formula above holds for any other $\mu_r$-variety. 

In the computations shown in this work, we focus on $\mu_2\cong\mathbb Z_2$. Let $1$ and $\varepsilon$ denote the trivial and nontrivial rational representations of $\mu_2$, respectively. Since $$ R_{\mathbb Q}(\mathbb Z_2)=\mathbb Z 1\oplus\mathbb Z \varepsilon,$$
we may write uniquely
$$[X]^{\mathbb Z_2}=[X]^+\otimes 1+[X]^-\otimes \varepsilon,$$
where $[X]^+=[X/\mathbb Z_2]$ and $[X]^-=[X]-[X/\mathbb Z_2]$. Indeed, the first coefficient is determined by taking $\mathbb Z_2$-invariants, while restricting to the trivial subgroup gives $[X]=[X]^++[X]^-$ (see \cite[Equation~(3.12)]{vogelthesis}).

The following lemma provides a product formula for equivariant virtual classes that will be useful in computations shown in Section \ref{sec:SL2}, following \cite[Equation~(3.12), Example~3.6.17 and Lemma~3.6.18]{vogelthesis}.

\begin{lema}\label{lem:plusminusformula}
    Let $X$ and $Y$ be algebraic varieties each endowed with an algebraic action of $\mathbb Z_2$. Assume that the multiplicativity hypothesis is satisfied, that is, $[X]\in V_{\mathbb Z_2}$ or $[Y]\in V_{\mathbb Z_2}$. Then
\begin{equation}
    [X\times Y]^+=[X]^+[Y]^+ + [X]^-[Y]^-, \quad [X\times Y]^-=[X]^+[Y]^- + [X]^-[Y]^+.
\end{equation}
\end{lema}

\begin{proof}
By the definition of $V_{\mathbb Z_2}$,
$$[X\times Y]^{\mathbb Z_2} = [X]^{\mathbb Z_2}[Y]^{\mathbb Z_2}.$$
Expanding both sides and using $\varepsilon^2=1$, the result follows by comparing the coefficients of $1$ and $\varepsilon$.
\end{proof}

The same construction and product formula apply to relative Grothendieck rings over a variety $S$ with trivial $\mathbb Z_2$-action, replacing Cartesian products by fiber products over $S$, as follows from \cite[Lemma~3.6.18 and Theorem~3.6.19]{vogelthesis}.

\subsection{The quantization functor}

The quantization functor, introduced in \cite{GONZALEZPRIETOlaxmonoTQFTrepvars}, is a lax monoidal functor
\[
\mathcal{Q}:\Span{\Var{\bk}} \longrightarrow \Mod{\K{\Var{\bk}}}
\]
which assigns to each variety $X$ its relative Grothendieck ring $\K{\Varrel{X}}$. For completeness, we recall the definitions needed for its construction.

Any morphism of algebraic varieties $f: S_1 \to S_2$ induces a functor $f^*: \Varrel{S_2} \to \Varrel{S_1}$ given by the fiber product along $f$. On the level of the Grothendieck rings, keeping the same symbols, this functor descends to a ring homomorphism $f^*: \K{\Varrel{S_2}} \to \K{\Varrel{S_1}}$ given by $f^*[X]_{S_2}=[X \times_{S_2} S_1]_{S_1}$. Considering the projection morphism $i:S \to \star$, the ring homomorphism $i^*:\K{\Var{\bk}} \to \K{\Varrel{S}}$ endows $\K{\Varrel{S}}$ with the structure of a $\K{\Var{\bk}}$-module.

Furthermore, $f$ induces a functor $f_!: \Varrel{S_1} \to \Varrel{S_2}$ given by composition with $f$, which descends to a $\K{\Var{\bk}}$-module morphism $f_![(X,h)]_{S_1}=[(X,f\circ h)]_{S_2}$ that we also denote by $f_!$. 

\begin{obs}
    We shall outline some important properties of the induced morphisms considered above. 
    \begin{itemize}
        \item Suppose we have a pullback of algebraic varieties
        \[
\begin{tikzcd}
S_1 \times_S S_2 \arrow[r, "g'"] \arrow[d, "f'"] & S_1 \arrow[d, "f"] \\
S_2 \arrow[r, "g"'] & S\end{tikzcd}
\]
Then, the so-called Beck–Chevalley property holds, that is, $g^* \circ f_! = (f')_! \circ (g')^*$.
\item Let $S$ be an algebraic variety with the decomposition given by $S= T \sqcup U$, where $i:T \hookrightarrow S$ is a closed embedding and $j:U \hookrightarrow S$ is an open embedding. Then, we have that $i_!i^*+j_!j^*: \K{\Varrel{S}} \to \K{\Varrel{S}}$ is the identity map.
    \end{itemize}
\end{obs}

We define the quantization functor by combining these two induced morphisms in the following way:

\begin{defi}
    The \emph{quantization functor} $\mathcal{Q}:\Span{\Var{\bk}} \rightarrow \Mod{\K{\Var{\bk}}}$ is defined on objects as $X \mapsto \K{\Varrel{X}}$ and on morphisms as
\begin{equation*}
\mathcal{Q}\left(X \xleftarrow{\ p\ } S \xrightarrow{\ q\ } Y\right)=\left(\K{\Varrel{X}}\xrightarrow{\ q_!p^*\ }\K{\Varrel{Y}}\right).
\end{equation*}
\end{defi}

The quantization functor is not strong monoidal but lax, since the natural map
$$
\Delta_{S_1,S_2}:
\K{\Varrel{S_1}}
\otimes_{\K{\Var{\bk}}}
\K{\Varrel{S_2}}
\longrightarrow
\K{\Varrel{S_1\times S_2}}
$$
is not an isomorphism in general. This map is induced by the external product of varieties, namely,
$$
\Delta_{S_1,S_2}
\left(
[X_1\xrightarrow{f_1}S_1]
\otimes
[X_2\xrightarrow{f_2}S_2]
\right)
=
[X_1\times X_2\xrightarrow{f_1\times f_2}S_1\times S_2].
$$

\subsection{Quantization of the virtual class of \texorpdfstring{$G$}{G}-representation varieties}

We now have all the ingredients needed to define the TQFT. In what follows, we carry out this construction in general dimension $n$. We mimic the functor construction done in \cite{gonzalezmotivictqft} but for horizontally convergent free boundary bordisms.

\begin{defi}
    The horizontally convergent TQFT \[\mathcal{Z}: \overline{\Bd}_n^\mathrm{free} \to \Mod{\K{\Var{\bk}}}\] is the symmetric lax monoidal functor defined as the composition $\mathcal{Q} \circ \mathcal{F}_G$.
    \label{def:tqft}
\end{defi}

For the sake of completeness, we outline the proof of the following analogue of \cite[Thm. 4.9]{gonzalezmotivictqft}. The construction is the same in $\overline{\Bd}_n^\mathrm{free}$. For a variety $S$, we denote by $\Unit{S}:=[\mathrm{id}:S\to S]_S$ the unit of $\K{\Varrel{S}}$.

\begin{tma}\label{tma:TQFTcomptesVirtClass}
Let $G$ be an affine algebraic group and $n\geq1$, and let
\[\mathcal{Z}: \overline{\Bd}_n^\mathrm{free}\to \Mod{\K{\Var{\bk}}}\]
be the TQFT introduced in Definition \ref{def:tqft}. Then $\mathcal{Z}(\varnothing)=\K{\Var{\bk}}$ and, for every bordism $W:\varnothing\to\varnothing$ in $\overline{\Bd}_n^\mathrm{free}$,
\[\mathcal{Z}(W)(\Unit{\star})=[R_G(W)]\in\K{\Var{\bk}},\]
where $\Unit{\star}$ is the unit of $\K{\Varrel{\star}}\cong\K{\Var{\bk}}$. Equivalently, $\mathcal{Z}(W)$ is multiplication by $[R_G(W)]$.
\end{tma}

\begin{proof}
    First we notice that for the empty set, $\mathcal{F}_G(\varnothing):=R_G(\varnothing)=\star$ is the singleton variety, hence $\mathcal{Z}(\varnothing)=\mathcal{Q} \circ \mathcal{F}_G(\varnothing)=\mathcal{Q}(\star) = \K{\Var{\bk}}$. A compact $n$-dimensional manifold can be described as a bordism with empty horizontal boundary $W: \varnothing \to \varnothing$. The associated span $\star \xleftarrow{c} R_G(W) \xrightarrow{c} \star$ induces a $\K{\Var{\bk}}$-module homomorphism $\mathcal{Z}(W)=c_!c^*:  \K{\Var{\bk}} \to  \K{\Var{\bk}}$. Using that $c^*$ is a ring homomorphism, we have that $c^*(\Unit{\star})=\Unit{R_G(W)}$ so that
     \begin{align*}
	\mathcal{Z}(W)(\Unit{\star}) & = c_!c^*(\Unit{\star}) = c_! \Unit{R_G(W)} \\
	&= c_![R_G(W)]_{R_G(W)} = [R_G(W) ]_\star = [R_G(W)].
     \end{align*}
This proves that the TQFT computes the virtual classes of $G$-representation varieties, as claimed.
\end{proof}

\begin{obs}\label{obs:no-basepoints}
     We note that the construction provided in \cite{gonzalezmotivictqft} uses fundamental groupoids and allows gluings involving several incoming or outgoing boundary components. Here, marked boundary points and connecting path classes are retained as part of the objects and morphisms of $\overline{\Bd}_n^\mathrm{free}$. Horizontal convergence allows each gluing to be treated by the ordinary fundamental group and Seifert--van Kampen theorem. Thus the apex uses one basepoint per connected component, and no additional factors of $[G]$ arise in the computation of $\mathcal{Z}(W)(\Unit{\star})$ for $W:\varnothing\to\varnothing$. Moreover, this formulation does not introduce additional subtleties such as conjugation maps in the corresponding spans; see \cite{GONZALEZPRIETOlaxmonoTQFTrepvars}.
\end{obs}

\section{Quantization of the Virtual Classes of Representation Varieties of Torus Links}\label{sec:QuantizationVirtualClasses}

In the previous section, we reviewed the construction of a TQFT that quantizes the virtual classes of representation varieties. We now specialize the previous construction to the representation varieties of torus links. In particular, we work with the category of horizontally convergent 3-bordisms with free boundary, marked points and connecting paths, $\overline{\Bd}_3^\mathrm{free}$. More precisely, we construct the relevant decomposition into bordisms of the complement of the torus link in $\sS^3$ and apply the lax symmetric monoidal functor \[\mathcal{Z}: \overline{\Bd}_3^\mathrm{free} \to \Mod{\K{\Var{\bk}}}\] defined in Definition \ref{def:tqft} to obtain the virtual class $[R_G(L^d_{n,m})]$. Then, thanks to the TQFT, we can exploit topological recursion to decompose the bordism associated with the torus link into simpler pieces, following the construction of Proposition \ref{prop:toruslinksarenested}, and compute the virtual class of its representation variety by composing the corresponding TQFT operators.

\subsection{Construction of the 3-bordism}\label{subsec:3bordism}

In this subsection, we aim to decompose the manifold $\sS^3-L^d_{n,m}$ into certain elementary pieces. Let $W$ be a compact bordism and $L \subseteq \operatorname{int}(W)$ a link. Choose a closed tubular neighborhood $\nu (L) \cong L \times D^2$, disjoint from the boundary (and free boundaries as well). Define its exterior as
$$E(W,L) := W - \operatorname{int}\nu(L).$$
This is a compact bordism with the same horizontal boundary as $W$ and with additional free boundary given by $\partial \nu (L) \cong \bigsqcup_{K\in \pi_0 (L)} \mathbb{T}$, where $\mathbb{T}=\sS^1\times \sS^1$ denotes the standard torus.

The radial deformation of each punctured normal disc onto its boundary extends to a strong deformation retraction of $W-L$ onto $E(W,L)$, fixing the horizontal boundary. This yields the following result:

\begin{lema}\label{lem:def_retract}
The inclusion $E(W,L)\hookrightarrow W-L$ induces an isomorphism of fundamental groups on each connected component and, consequently, an isomorphism of representation varieties $$R_G(E(W,L))\cong R_G(W-L),$$ compatible with the restriction morphisms induced by the horizontal boundary inclusions.
\end{lema}

\begin{figure}[htbp]
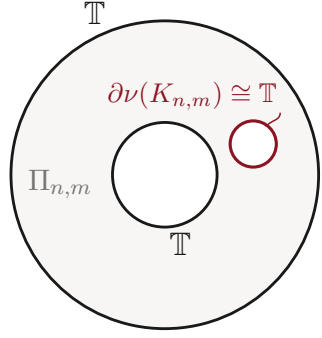

    \centering
    \vspace{-1.5ex}
    \resizebox{0.28\linewidth}{!}{\thickenedtoruscut}
    \vspace{-2ex}
    \caption{Transverse section of the $(n,m)$-fat torus associated with the $(n,m)$-torus knot.}
    \label{fig:thickened-torus-cut}
\end{figure}

Recall the decomposition of $\sS^3$ into two solid tori, described using hyperspherical coordinates in Section \ref{sec:rep_varieties}. Using this decomposition, we now introduce the relevant bordisms in the category $\overline{\Bd}_3^\mathrm{free}$. They are the following (recall that $\mathbb{T}=\sS^1\times \sS^1$ denotes the standard torus):

\begin{itemize}
    \item The \emph{outer solid torus}, denoted by $\mathbf{T}:\varnothing\to \mathbb{T}$, which corresponds to one of the solid tori of the aforementioned decomposition.
    \item The \emph{$(n,m)$-fat torus}, denoted by $\Pi_{n,m}:\mathbb{T} \to \mathbb{T}$. It is the exterior of a tubular neighborhood of one copy of $K_{n,m}$ in $\mathbb{T}\times[0,1]$, with $\partial \nu(K_{n,m})$ as a free boundary. A diagram of the transverse section of $\Pi_{n,m}$ can be seen in Figure \ref{fig:thickened-torus-cut}. 
    \item The \emph{inner solid torus}, denoted by $\mathbf{T}^\dagger:\mathbb{T}\to\varnothing$, which corresponds to the other solid torus of the decomposition.
\end{itemize}

We fix a marked point on the torus $\mathbb{T}$ and use its images under the boundary parametrizations. On $\Pi_{n,m}$ we fix the connecting path used in the Seifert--van Kampen computation of Proposition \ref{pi1tori}. Repeated copies carry the same data, and their paths are concatenated under gluing.

\subsection{Computations of the operators}
Thanks to Lemma \ref{lem:def_retract} and Proposition \ref{prop:toruslinksarenested}, the exterior $E(\sS^3,L_{n,m}^d)$ of the torus link is boundary-preserving diffeomorphic to the composition of bordisms
$$E(\sS^3,L_{n,m}^d)=\mathbf{T}^{\dagger}\circ\Pi_{n,m}^d\circ\mathbf{T},$$
where $\mathbf{T}:\varnothing\to \mathbb{T}$ is the outer solid torus, $\Pi_{n,m}:\mathbb{T}\to\mathbb{T}$ is the $(n,m)$-fat torus and $\mathbf{T}^\dagger:\mathbb{T}\to\varnothing$ is the inner solid torus. Hence, to compute the virtual class $[R_G(L_{n,m}^d)]$ of the representation variety of the torus link, we only need to compute the operators
\begin{align*}
\mathcal{Z}(\mathbf{T}) &: \K{\Var{\bk}} \longrightarrow \K{\Varrel{\comvar}},\\
\mathcal{Z}(\mathbf{T}^\dagger) &: \K{\Varrel{\comvar}} \longrightarrow \K{\Var{\bk}},\\
\mathcal{Z}(\Pi_{n,m}) &: \K{\Varrel{\comvar}} \longrightarrow \K{\Varrel{\comvar}},
\end{align*}
where
$$\comvar:=R_G(\mathbb{T})=\{(g_1,g_2)\in G^2 \ | \ g_1g_2=g_2g_1\}.$$

\begin{figure}[ht]
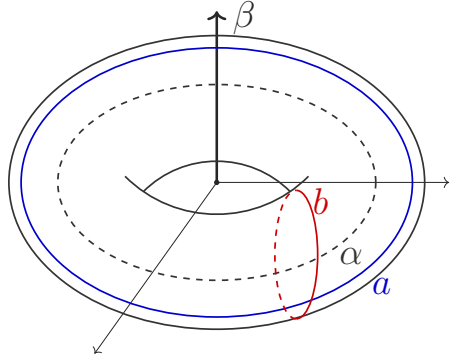

\centering
\resizebox{0.35\linewidth}{!}{\torusinsphere}
\caption{Generators of $\pi_1(\mathbb{T}),\pi_1(\mathbf{T})$ and $\pi_1(\mathbf{T}^\dagger)$.}
    \label{fig:fund_group_of_torus}
\end{figure}

Let us start with the outer solid torus $\mathbf{T}$. The fundamental groups of the manifolds appearing in the bordism are $\pi_1(\mathbf{T})=\langle \beta\rangle$ and $\pi_1(\mathbb{T})=\langle a,b\ | \ [a,b]\rangle$, where $a$ is a meridian of the torus and $b$ is a parallel. This gives us the representation varieties $R_G(\varnothing)=\star, R_G(\mathbf{T})=G$ and $R_G(\mathbb{T})=\comvar$. Consider the inclusion map $i:\mathbb{T}\hookrightarrow\mathbf{T}$. It is easy to see that $i(a)$ is contractible in $\mathbf{T}$, whereas $i(b)$ is homotopic to $\beta$ (see Figure \ref{fig:fund_group_of_torus}). Therefore, the map induced on fundamental groups $i_*:\pi_1(\mathbb{T})\to \pi_1(\mathbf{T})$ is given by $i_*(a)=e$ and $i_*(b)=\beta$, and, hence, the induced map on representation varieties is just $q_2:G\to\comvar, g\mapsto (e,g)$. Hence, the geometrization step of the TQFT assigns to the bordism $\mathbf{T}:\varnothing \to \mathbb{T}$ the span
\begin{center}
\begin{tikzcd}[ampersand replacement=\&]
\& G \arrow[dl, "c"'] \arrow[dr, "q_2"] \& \\
\star \& \& \comvar
\end{tikzcd}
\end{center}
where $c:G\to \star$ is the constant map and $q_2:G\to\comvar, g\mapsto (e,g)$. Now, the quantization gives us
\begin{equation}\label{eq:op_Z_T}
\mathcal{Z}(\mathbf{T})=(q_2)_!c^*:\K{\Var{\bk}}\to \K{\Varrel{\comvar}}.
\end{equation}

A similar argument gives us the operator of the inner solid torus
\begin{equation}\label{eq:op_Z_T_dagger}
\mathcal{Z}(\mathbf{T}^\dagger)=c_!q_1^*:\K{\Varrel{\comvar}}\to\K{\Var{\bk}},
\end{equation}
where $q_1:G\to\comvar$ is given by $g\mapsto(g,e)$, as, in $\mathbf{T}^\dagger$, the loop $a$ is homotopic to the generator $\alpha$ of $\pi_1(\mathbf{T}^\dagger)$ and $b$ is nullhomotopic.

We finally address the case of the $(n,m)$-fat torus. For this, we will first compute the fundamental group of this piece:

\begin{prop}\label{pi1tori}
    The fundamental group of the $(n,m)$-fat torus is 
    $$\pi_1(\Pi_{n,m})=\langle a_1,b_1,a_2,b_2 \ | \ [a_1,b_1], \ [a_2,b_2], \ a_1^nb_1^ma_2^{-n}b_2^{-m}\rangle.$$ Moreover, the inclusions of the incoming and outgoing horizontal boundaries induce the homomorphisms $$ \begin{aligned}
(\iota_-)_*:\pi_1(\mathbb{T})&\longrightarrow\pi_1(\Pi_{n,m}),
& (a,b)&\longmapsto (a_1,b_1),\\
(\iota_+)_*:\pi_1(\mathbb{T})&\longrightarrow\pi_1(\Pi_{n,m}),
& (a,b)&\longmapsto (a_2,b_2),
\end{aligned}
$$ where $a,b$ are the meridian and parallel generators of $\pi_1(\mathbb{T})$, respectively.
\end{prop}

\begin{figure}[ht]
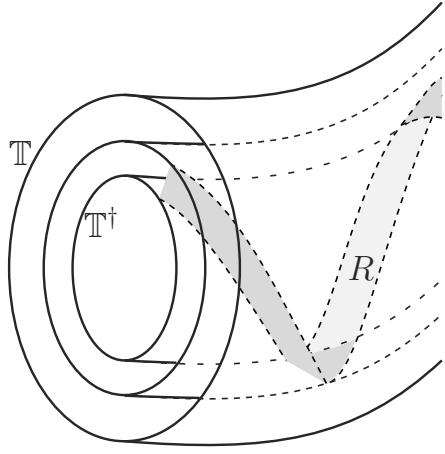

\centering
\resizebox{0.35\linewidth}{!}{\vankampen}
\caption{Schematic representation of the $(n,m)$-fat torus $\Pi_{n,m}$ used in the Seifert--van Kampen decomposition.}
    \label{fig:van_kampen}
\end{figure}

\begin{proof}
    We will use Seifert--van Kampen's theorem to prove the result. By Lemma \ref{lem:def_retract}, we may compute the fundamental group using the knot-complement $P_{n,m}:=(\mathbb{T}\times[0,1])-K_{n,m}$, of which $\Pi_{n,m}$ is a deformation retract. Let us denote by $\mathbb{T}$ the external torus of $\Pi_{n,m}$ and $\mathbb{T}^{\dagger}$ the inner one, so that $\partial^-_{\mathrm h}\Pi_{n,m}\sqcup \partial^+_{\mathrm h}\Pi_{n,m}=\mathbb{T}\sqcup\mathbb{T}^\dagger$. Let also $\mathbb{T}_{n,m}$ be the torus between $\mathbb{T}$ and $\mathbb{T}^\dagger$ (and concentric with them) that contains the $(n,m)$-torus knot $K_{n,m}$ that we remove to obtain $P_{n,m}$. Let $\varepsilon$ be the distance between $\mathbb{T}^{\dagger}$ and $\mathbb{T}_{n,m}$. Consider the normal vector field $\nu$ of length $\varepsilon$ along the torus knot $K_{n,m}$ that is also normal to $\mathbb{T}_{n,m}$ and points inwards. In other words, $\nu$ assigns to each point of $K_{n,m}$ the vector normal to $\mathbb{T}_{n,m}$ that starts at that point and ends at a point of $\mathbb{T}^{\dagger}$. The image of $\nu$ is a copy of the torus knot $K_{n,m}$ but in the inner torus $\mathbb{T}^\dagger$, and the union of all the normal vectors is a ribbon $R$ whose boundary components are $K_{n,m}$ and the mentioned copy in $\mathbb{T}^\dagger$ (see Figure \ref{fig:van_kampen}). Let $A=P_{n,m}-(\mathbb{T}^\dagger\cup R)$ be the open subset of $P_{n,m}$ obtained by removing the inner torus $\mathbb{T}^\dagger$ and the ribbon $R$. The open set $A$ retracts to the outer torus $\mathbb{T}$, so its fundamental group is 
    $$\pi_1(A)=\pi_1(\mathbb{T})=\langle a_1,b_1 \ | \ [a_1,b_1]\rangle,$$
    where $a_1$ and $b_1$ are the meridian and parallel of $\mathbb{T}$ respectively, in the boundary coordinates fixed above.

    Analogously, we can define the open set $B=P_{n,m}-(\mathbb{T}\cup S)$, where now $S$ is a ribbon whose boundary components are $K_{n,m}$ and a copy of it in $\mathbb{T}$. Now, $\pi_1(B)=\langle a_2,b_2 \ | \ [a_2,b_2]\rangle$ is generated by a meridian and parallel of $\mathbb{T}^\dagger$.

    Now, the intersection $U=A\cap B$ of the two open sets is just the interior of $P_{n,m}$ with the two ribbons $R$ and $S$ removed, i.e.\ $U=\mathring{P}_{n,m}-(R\cup S)$. This subset retracts to $\mathbb{T}_{n,m}-K_{n,m}$, so its fundamental group is
    $$\pi_1(U)=\pi_1(\mathbb{T}_{n,m}-K_{n,m})$$
    which is generated by a copy $K'_{n,m}$ of the torus knot $K_{n,m}$, seen as a loop in the fundamental group.

    Finally, Seifert--van Kampen theorem shows that the fundamental group of the $(n,m)$-fat torus is 
    $$\pi_1(\Pi_{n,m})=\langle a_1,b_1,a_2,b_2 \ | \ [a_1,b_1], \ [a_2,b_2], \ (i_A)_*(K'_{n,m})(i_B)_*(K'_{n,m})^{-1}\rangle,$$
    where $i_A:U\hookrightarrow A$ and $i_B:U\hookrightarrow B$ are the inclusion maps and $(i_A)_*$ and $(i_B)_*$ are the induced ones at the level of fundamental groups. But $K'_{n,m}$ is an $(n,m)$-torus knot, so we have $(i_A)_*(K'_{n,m})=a_1^nb_1^m$ and $(i_B)_*(K'_{n,m})=a_2^nb_2^m$, giving the stated presentation. Moreover, the boundary parametrizations identify the generators $a,b$ of $\pi_1(\mathbb{T})$ with $a_1,b_1$ on the exterior torus and with $a_2,b_2$ on the inner torus, respectively. This gives the stated homomorphisms $(\iota_-)_*$ and $(\iota_+)_*$.
\end{proof}

\begin{coro}
    The representation variety of the $(n,m)$-fat torus is 
    $$R_G(\Pi_{n,m})=\{(g_1,h_1,g_2,h_2)\in G^4 \ | \ [g_1,h_1]= [g_2,h_2]=e, \ g_1^nh_1^m = g_2^nh_2^m\}.$$ Moreover, the restriction morphisms induced by the horizontal boundary inclusions are $$
    \begin{aligned}
    p_1:R_G(\Pi_{n,m})&\longrightarrow\comvar,
    & (g_1,h_1,g_2,h_2)&\longmapsto(g_1,h_1),\\
    p_2:R_G(\Pi_{n,m})&\longrightarrow\comvar,
    & (g_1,h_1,g_2,h_2)&\longmapsto(g_2,h_2).
    \end{aligned}
    $$
\end{coro}

Thus, the geometrization functor assigns to $\Pi_{n,m}$ the span
\begin{center}
\begin{tikzcd}[ampersand replacement=\&]
\& R_G(\Pi_{n,m})
\arrow[dl, "p_1"']
\arrow[dr, "p_2"] \& \\
\comvar \& \& \comvar.
\end{tikzcd}
\end{center}

Finally, we apply the quantization functor to obtain
\begin{equation}\label{eq:Z_tori}
    \mathcal{Z}(\Pi_{n,m})=(p_2)_!(p_1)^*:\K\Varrel{ \comvar}\to\K\Varrel{\comvar}.
\end{equation}

\subsection{Computation of the virtual class}\label{sec:computation_virtualclass} We are now in a position to apply the TQFT to study the virtual class of the representation variety $R_G(L^d_{n,m})$ of a torus link. By Lemma \ref{lem:def_retract} and Theorem \ref{tma:TQFTcomptesVirtClass} we have
$$[R_G(L^d_{n,m})]=\mathcal{Z}(E(\sS^3,L_{n,m}^d))(\Unit{\star}),$$
where $\Unit{\star}=[\star]\in\K\Var{\bk}$ is the identity element. By the decomposition of $E(\sS^3,L^d_{n,m})$ into bordisms, we have
$$\mathcal{Z}(E(\sS^3,L^d_{n,m}))(\Unit{\star})=\mathcal{Z}(\mathbf{T}^\dagger)\circ\mathcal{Z}(\Pi_{n,m})^d\circ\mathcal{Z}(\mathbf{T})(\Unit{\star}).$$

The first computation is easy, since by equation (\ref{eq:op_Z_T}), we have 
$$\mathcal{Z}(\mathbf{T})(\Unit{\star})=(q_2)_!c^*(\Unit{\star})=(q_2)_!(\Unit{G})=(q_2)_!([\id_G:G\to G]_G)=[q_2:G\to \comvar]_{\comvar}.$$
Here $c^*(\Unit{\star})=\Unit{G}\in\K\Varrel{G}$ because $c^*$ is a ring homomorphism.

Let us now consider $\mathcal{Z}(\Pi_{n,m})$. Equation (\ref{eq:Z_tori}) gives us
$$\mathcal{Z}(\Pi_{n,m})([q_2:G\to\comvar]_{\comvar})=(p_2)_!(p_1)^*([q_2:G\to\comvar]_{\comvar}).$$
The image of $[q_2:G\to\comvar]_{\comvar}$ under $(p_1)^*$ is the class in $\K\Varrel{R_G(\Pi_{n,m})}$ of the upper horizontal morphism in the pullback diagram
\begin{center}
\begin{tikzcd}[ampersand replacement=\&]
X\arrow[r] \arrow[d] \& R_G(\Pi_{n,m})\arrow[d, "p_1"]\\
G\arrow[r, "q_2"]  \& \comvar
\end{tikzcd}
\end{center}
where
\begin{equation}\label{eq:X}
X:=G\times_{\comvar}R_G(\Pi_{n,m})=\{(h_1,g_2,h_2)\in G^3 \ | \ [g_2,h_2]=e, h_1^m=g_2^{n}h_2^{m}\}.
\end{equation}
Applying $(p_2)_!$ to this morphism gives us
$$(p_2)_!([X\to R_G(\Pi_{n,m})]_{R_G(\Pi_{n,m})})=[\varphi:X\to \comvar]_{\comvar}$$
where the morphism $\varphi$ is given by $(h_1,g_2,h_2)\mapsto (g_2,h_2)$.

The class of $\varphi:X\to \comvar$ is not known a priori, and applying the operator $\mathcal{Z}(\Pi_{n,m})$ to it only leads to more complicated expressions, so this does not facilitate the computation of $\mathcal{Z}(\Pi_{n,m})^d$. To overcome this difficulty, we are going to stratify the variety $X$ in terms of the stabilizers of $h_1^m$ and $h_1$ under conjugation by $G$. The reason for choosing this stratification will become clear in a moment. Let us first consider the case where $G$ is a connected reductive algebraic group. The quotients used below are geometric quotients of free actions, denoted by $/$, rather than affine GIT quotients. With minor modifications, the arguments extend to arbitrary affine algebraic groups, as discussed in Remark \ref{rem:general-algebraic-groups}.

Let $H_1 \leq H_2 \leq G$ be two closed subgroups of $G$. Consider the locus
$$X_{H_1,H_2}:=\{(h_1,g_2,h_2)\in X \ | \ (\Stab_G(h_1),\Stab_G(h_1^m))\sim (H_1,H_2)\}$$
where stabilizers are taken for conjugation and $\sim$ denotes simultaneous conjugation of the two subgroups. $G$ is reductive, so it has a finite partition into Jordan classes, and each Jordan class is a smooth locally closed subvariety. Thus, $X_{H_1,H_2}$ is a finite union of locally closed subvarieties. Let
$$
\mathrm{PSb}(G):=\{(H_1,H_2)\mid H_1\leq H_2\leq G,\,X_{H_1,H_2}\neq\varnothing\},
\qquad
\mathcal{PS}b(G):=\mathrm{PSb}(G)/\sim.
$$

In this setting, $X$ admits the following stratification:
\begin{equation}
X=\bigsqcup_{[(H_1,H_2)]\in\mathcal{PS}b(G)}X_{H_1,H_2}.
\label{eq:stratification}
\end{equation}

Indeed, since $G$ is reductive, only finitely many simultaneous conjugacy classes of stabilizer pairs $\bigl(\Stab_G(h),\Stab_G(h^m)\bigr)$ occur in the stratification. To see this, let $h=su$ be the Jordan decomposition of $h$. Then $\Stab_G(h)=\Stab_G(s)\cap\Stab_G(u)$ and, since $u$ and $u^m$ have the same centralizer in characteristic zero, $\Stab_G(h^m)=\Stab_G(s^m)\cap\Stab_G(u)$. Up to conjugation, $s$ lies in a fixed maximal torus, and then the pair $\bigl(\Stab_G(s),\Stab_G(s^m)\bigr)$ takes only finitely many values \cite[Theorem~2.2]{humphreys1995conjugacy}. Finally, $u$ is a unipotent element of the reductive group $\Stab_G(s)$, which has finitely many conjugacy classes of unipotent elements \cite{humphreys1995conjugacy}, and conjugation by $\Stab_G(s)\subseteq\Stab_G(s^m)$ conjugates both stabilizers simultaneously. By the additivity of relative virtual classes applied to the stratification \eqref{eq:stratification}, we have
$$
[\varphi:X\to\comvar]_{\comvar}
=
\sum_{[(H_1,H_2)]}
[\varphi_{H_1,H_2}:X_{H_1,H_2}\to\comvar]_{\comvar},
$$
where $\varphi_{H_1,H_2}$ denotes the restriction of $\varphi$ to $X_{H_1,H_2}$.

We would like the submodule $$\mathcal{W}_0=\langle[X_{H_1,H_2}\to\comvar]_{\comvar}\rangle_{[(H_1,H_2)]\in\mathcal{PS}b(G)}$$ 
to be invariant under $\mathcal{Z}(\Pi_{n,m})$, since, in that case, the image of $\mathcal{Z}(\Pi_{n,m})^{d-1}([\varphi :X\to \comvar]_{\comvar})$ for any $d$ would be completely determined by the images $\mathcal{Z}(\Pi_{n,m})([X_{H_1,H_2}\to \comvar]_{\comvar})$. Unfortunately, this is in general not the case. For this reason, we also need to introduce
$$X_{H_1,H_2}^0:=\{(h_1,g_2,h_2)\in X \ | \ (\Stab_G(h_1),\Stab_G(h_1^m))= (H_1,H_2)\}$$
and the morphisms
$$\varphi_{H_1,H_2}^0:X_{H_1,H_2}^0\times G\to \comvar, (h_1,g_2,h_2,g)\mapsto (gg_2g^{-1},gh_2g^{-1}).$$

To relate $X^0_{H_1,H_2}$ to the original strata $X_{H_1,H_2}$, we consider the quotients of $X^0_{H_1,H_2}\times G$ by $H_1$ and by the subgroup preserving the pair $(H_1,H_2)$.

Let $H_1\leq H_2\leq G$, and set $M_{H_1,H_2}=N_G(H_1)\cap N_G(H_2)$ and $N_{H_1,H_2}=M_{H_1,H_2}/H_1$, where $N_G(H)=\{g\in G\mid gHg^{-1}=H\}$ denotes the normalizer of $H$ in $G$. We write 
$$
\widehat X^0_{H_1,H_2}:=(X^0_{H_1,H_2}\times G)/M_{H_1,H_2},
\qquad
\widetilde X^0_{H_1,H_2}:=(X^0_{H_1,H_2}\times G)/H_1,
$$
where the group $M_{H_1,H_2}$ acts on $X^0_{H_1,H_2}\times G$ by
$$
r\cdot(h_1,g_2,h_2,g)
=
(rh_1r^{-1},rg_2r^{-1},rh_2r^{-1},gr^{-1}),
\qquad r\in M_{H_1,H_2}.
$$

The quotient $\widehat X^0_{H_1,H_2}$ provides a model for the original stratum $X_{H_1,H_2}$, while $\widetilde X^0_{H_1,H_2}$ is an auxiliary variety introduced to simplify the study of the TQFT operator and relate its action to the virtual classes of the original strata.

The action of $H_1$ used to define $\widetilde X^0_{H_1,H_2}$ is the restriction of this action. Since $H_1\trianglelefteq M_{H_1,H_2}$ is a normal subgroup, the $M_{H_1,H_2}$-action descends to an action of $N_{H_1,H_2}=M_{H_1,H_2}/H_1$ on $\widetilde X^0_{H_1,H_2}$, given by $$
rH_1\cdot[h_1,g_2,h_2,g]=[rh_1r^{-1},rg_2r^{-1},rh_2r^{-1},gr^{-1}].
$$

\begin{prop}\label{prop:quotient-torsor}
Let $M=M_{H_1,H_2}$. For $K=H_1$ or $M$, consider
$$
Y_K=\{(gK,x)\in(G/K)\times G^3 \mid g^{-1}\cdot x\in X^0_{H_1,H_2}\}.
$$
Then the geometric quotients admit the explicit descriptions
$$
\widetilde X^0_{H_1,H_2}\cong Y_{H_1}, \qquad \widehat X^0_{H_1,H_2}\cong Y_M.
$$
Moreover, the natural projection $\widetilde X^0_{H_1,H_2}\to\widehat X^0_{H_1,H_2}$ is an $N_{H_1,H_2}$-torsor, and we have $\widehat X^0_{H_1,H_2} =\widetilde X^0_{H_1,H_2}/N_{H_1,H_2}$.
\end{prop}

\begin{proof}
Since $X^0_{H_1,H_2}$ is constructible, choose a finite $M$-invariant locally closed decomposition and work on
each piece.

First, consider $K=H_1$. The definition of $Y_{H_1}$ does not depend on the representative $g$ of $gH_1$, since $X^0_{H_1,H_2}$ is invariant under conjugation by $H_1$. The map
$$
X^0_{H_1,H_2}\times G\longrightarrow Y_{H_1}, \qquad (x,g)\longmapsto(gH_1,g\cdot x),
$$
is $H_1$-invariant and identifies $Y_{H_1}$ with the geometric quotient $(X^0_{H_1,H_2}\times G)/H_1$. The same argument applies to $K=M$, since $X^0_{H_1,H_2}$ is also $M$-invariant.

Since $H_1\trianglelefteq M$ is a normal subgroup, the projection $G/H_1\to G/M$ is an $M/H_1=N_{H_1,H_2}$-torsor. Moreover, the natural projection $Y_{H_1}\to Y_M$, given by $(gH_1,x)\mapsto(gM,x)$, fits into the fiber-product identification
$$
Y_{H_1}\cong Y_M\times_{G/M}(G/H_1).
$$
Thus, $Y_{H_1}\to Y_M$ is the base change of the $N_{H_1,H_2}$-torsor $G/H_1\to G/M$, and is therefore itself an $N_{H_1,H_2}$-torsor. Consequently, $\widehat X^0_{H_1,H_2} =\widetilde X^0_{H_1,H_2}/N_{H_1,H_2}$.
\end{proof}

\begin{tma}\label{tma:covering}
Let $H_1\leq H_2\leq G$. There is a natural regular bijection
$$
\Phi_{H_1,H_2}:\widehat X^0_{H_1,H_2}
\longrightarrow X_{H_1,H_2},
$$
induced by simultaneous conjugation, namely
$$
\Phi_{H_1,H_2}[h_1,g_2,h_2,g]
=(gh_1g^{-1},gg_2g^{-1},gh_2g^{-1}).
$$
Moreover, the morphism $\varphi^0_{H_1,H_2}$ descends to morphisms $\widetilde{\varphi}^0_{H_1,H_2}$ and $\widehat{\varphi}^0_{H_1,H_2}$, making the following diagram commute:
\[
\begin{tikzcd}[column sep=large, row sep=large]
\widetilde X^0_{H_1,H_2}
    \arrow[r, "q"]
    \arrow[dr, "\widetilde{\varphi}^0_{H_1,H_2}"']
&
\widehat X^0_{H_1,H_2}
    \arrow[r, "\Phi_{H_1,H_2}"]
    \arrow[d, "\widehat{\varphi}^0_{H_1,H_2}"]
&
X_{H_1,H_2}
    \arrow[dl, "\varphi_{H_1,H_2}"]
\\
& \comvar &
\end{tikzcd}
\]
where $q$ is the $N_{H_1,H_2}$-torsor of Proposition \ref{prop:quotient-torsor}. In particular, $\Phi_{H_1,H_2}$ induces an equality of the corresponding relative virtual classes over $\comvar$.
\end{tma}

\begin{proof}
The morphism
\[
    X^0_{H_1,H_2}\times G\longrightarrow X_{H_1,H_2},\qquad
    (h_1,g_2,h_2,g)\longmapsto
    (gh_1g^{-1},gg_2g^{-1},gh_2g^{-1})
\]
is $M_{H_1,H_2}$-invariant and surjective, so it descends to a surjective morphism
\[
    \Phi_{H_1,H_2}:
    \widehat X^0_{H_1,H_2}\longrightarrow X_{H_1,H_2}.
\]

If two classes $[h_1,g_2,h_2,g]$ and $[h_1',g_2',h_2',g']$ have the same image, then $r=(g')^{-1}g$ satisfies $(h_1',g_2',h_2')=(rh_1r^{-1},rg_2r^{-1},rh_2r^{-1})$. Since both triples lie in $X^0_{H_1,H_2}$, we get $rH_1r^{-1}=\Stab_G(h_1')=H_1$ and $rH_2r^{-1}=\Stab_G((h_1')^m)=H_2$. Hence
\[
    r=(g')^{-1}g\in N_G(H_1)\cap N_G(H_2)=M_{H_1,H_2},
\]
and $r\cdot(h_1,g_2,h_2,g)=(h_1',g_2',h_2',gr^{-1})=(h_1',g_2',h_2',g')$, so the two classes coincide. Thus $\Phi_{H_1,H_2}$ is bijective. 
    
Finally, the morphism $\varphi^0_{H_1,H_2}$ is invariant under the relevant actions and hence descends to $\widetilde{\varphi}^0_{H_1,H_2}$ and $\widehat{\varphi}^0_{H_1,H_2}$. The commutativity of the diagram follows from
$$
\widetilde{\varphi}^0_{H_1,H_2}
=\widehat{\varphi}^0_{H_1,H_2}\circ q,
\qquad
\widehat{\varphi}^0_{H_1,H_2}
=\varphi_{H_1,H_2}\circ\Phi_{H_1,H_2},
$$
since all three morphisms are induced by $(h_1,g_2,h_2,g)\mapsto (gg_2g^{-1},gh_2g^{-1})$. The equality of relative virtual classes is a consequence of Proposition \ref{prop:equal_virtual_class}.
\end{proof}

As we will see in Theorem \ref{tma:resultados} and Corollary \ref{coro:associated-centralizer-eigenvector}, under suitable hypotheses, the classes of the morphisms $\widetilde{\varphi}^0_{H_1,H_2}:\widetilde{X}^0_{H_1,H_2}\longrightarrow \comvar$ induced by the morphisms $\varphi_{H_1,H_2}^0:X^0_{H_1,H_2}\times G\to \comvar$ are eigenvectors of $\mathcal{Z}(\Pi_{n,m})$. Hence, being able to relate our original strata $X_{H_1,H_2}\to \comvar$ with these new varieties can be very useful for doing the computations of the TQFT. For example, suppose that $N_{H_1,H_2}$ is a \emph{special algebraic group}, namely, an algebraic group
whose torsors for the étale topology are Zariski locally trivial.

\begin{coro}
If $N_{H_1,H_2}$ is special, then,
\begin{equation*}
[\widetilde\varphi^0_{H_1,H_2}:\widetilde X^0_{H_1,H_2}\to\comvar]_{\comvar} =[N_{H_1,H_2}]
[\varphi_{H_1,H_2}:X_{H_1,H_2}\to\comvar]_{\comvar}.
\end{equation*}
\end{coro}

\begin{proof}
The $N_{H_1,H_2}$-torsor $\widetilde X^0_{H_1,H_2}\to\widehat X^0_{H_1,H_2}$ is Zariski locally trivial. Therefore, its relative virtual class factors as $[N_{H_1,H_2}]$ times the relative virtual class of $\widehat X^0_{H_1,H_2}\to\comvar$. Then, apply Theorem \ref{tma:covering} to identify the latter class with $[\varphi_{H_1,H_2}]$.
\end{proof}

However, the groups $N_{H_1,H_2}$ need not be special; for instance, a nontrivial finite group is never special. When the mentioned normalizer quotient $N_{H_1,H_2}$ is finite, as occurs in Section \ref{sec:SL2}, its action must still be retained in the quotient computations. In these cases, the use of \emph{equivariant virtual classes} allows us to still do the computations. 

Now we focus on proving the announced result that the classes of the morphisms $\widetilde{\varphi}_{H_1,H_2}^0:(X^0_{H_1,H_2}\times G) / H_1\to \comvar$ are, under suitable hypotheses, eigenvectors of $\mathcal{Z}(\Pi_{n,m})$, as well as a stronger result when $H_2=G$. To do this, we need two technical lemmata .

\begin{lema}\label{g2nh2m=stabh1m}
Let $h_1\in G$, and let
$$F(h_1):=\{(g_2,h_2)\in G^2\ | \ [g_2,h_2]=e, \ g_2^nh_2^m=h_1^m\}.$$
If $a,b\in\Z$ satisfy $an+bm=1$, the map
$$\Phi_{h_1}:\Stab_G(h_1^m)\to F(h_1), h\mapsto (h^mh_1^{am},h^{-n}h_1^{bm})$$
is an isomorphism.
\end{lema}

\begin{proof}
We first see that $\Phi_{h_1}$ is well defined. Indeed, we need to see that $h^mh_1^{am}$ and $h^{-n}h_1^{bm}$ commute and that $(h^mh_1^{am})^n(h^{-n}h_1^{bm})^m=h_1^m$. The commutativity is straightforward, since $h\in\Stab_G(h_1^m)$. The other equality also follows from the commutativity of $h$ with $h_1^m$:
$$(h^mh_1^{am})^n(h^{-n}h_1^{bm})^m=h^{mn}h_1^{amn}h^{-nm}h_1^{bmm}=h^{nm}h^{-nm}h_1^{m(an+bm)}=h_1^m.$$

Now, the inverse is the regular morphism $$\Phi^{-1}_{h_1}:F(h_1)\to \Stab_G(h_1^m), (g_2,h_2)\mapsto g_2^bh_2^{-a}.$$ The identities $an+bm=1$ and $g_2h_2=h_2g_2$ verify both compositions.
\end{proof}

\begin{lema}\label{lma:pullback_commutes_with_quot}
Let $q:P\to Q$ be an $H$-torsor for an algebraic group $H$, and let $f:Q\to S$ and $g:Y\to S$ be morphisms. Then $Q\times_S Y$ is the geometric quotient of $P\times_S Y$ by the induced action of $H$ on the first factor. In particular, the quotient exists and there is a canonical isomorphism
$$
(P\times_S Y)/H\cong Q\times_S Y.
$$
Moreover, for any $H$-invariant morphism $\varphi:P\times_S Y\to B$, this is an isomorphism over $B$, where the morphism $Q\times_S Y\to B$ is induced by $\varphi$.
\end{lema}

\begin{proof}
There is a canonical identification $ P\times_S Y \cong P\times_Q(Q\times_S Y)$, where $Q\times_S Y\to Q$ is the first projection. Hence $P\times_S Y\longrightarrow Q\times_S Y$ is the base change of the $H$-torsor $q:P\to Q$, and is therefore again an $H$-torsor. Consequently, $Q\times_S Y$ is its geometric quotient by $H$, giving the claimed isomorphism.

Since $\varphi$ is $H$-invariant, it descends to a unique morphism $\overline{\varphi}:Q\times_S Y\to B$. Thus, the quotient identification is compatible with the morphisms to $B$, and the corresponding relative virtual classes coincide.
\end{proof}

We now establish the action of the TQFT operator on the relative virtual classes associated with the stabilizer stratification.

\begin{tma}\label{tma:resultados} 

\begin{itemize}

\item[(a)] Set
$$
X_G:=\{(h_1,g_2,h_2)\in X
\mid\Stab_G(h_1^m)=G\},
\qquad \varphi_G=\varphi|_{X_G}.
$$
Then
$$
\mathcal Z(\Pi_{n,m})([\varphi_G])
=[G][\varphi_G].
$$

\item[(b)] Let $H_1\leq H_2\leq G$. Consider the morphism $\widetilde\varphi^0_{H_1,H_2}: \widetilde X^0_{H_1,H_2}\to\comvar$ of Theorem \ref{tma:covering}. Define the associated bundle
$$
\mathcal E_{H_1,H_2}
:=(X^0_{H_1,H_2}\times G\times H_2)/H_1
\xrightarrow{\ \pi\ }\widetilde X^0_{H_1,H_2},
$$
where $H_1$ acts by
$$
h\cdot(x,g,h')=(h\cdot x,gh^{-1},hh'h^{-1}),
$$
with $x=(h_1,g_2,h_2)\in X^0_{H_1,H_2}$,
$h\in H_1$, and $h\cdot x$ denoting simultaneous conjugation. Then, the operator is given by
\begin{equation}\label{eq:associated-centralizer-operator}
\mathcal Z(\Pi_{n,m})
([\widetilde\varphi^0_{H_1,H_2}])
=
[\widetilde\varphi^0_{H_1,H_2}\circ\pi:
\mathcal E_{H_1,H_2}\to\comvar]_{\comvar}.
\end{equation}
\end{itemize}
\end{tma}

\begin{proof}
For part (a), the pullback is
\begin{align*}
X_G\times_{\comvar}R_G(\Pi_{n,m})
=\{&(h_1,g_2,h_2,g_3,h_3)\in G^5\mid\\
&[g_2,h_2]=[g_3,h_3]=e,\\
&h_1^m=g_2^nh_2^m=g_3^nh_3^m,\ 
\Stab_G(h_1^m)=G\}.
\end{align*}
By Lemma \ref{g2nh2m=stabh1m}, the map
$$
(h_1,g_2,h_2,g_3,h_3)
\longmapsto
(h_1,g_3,h_3,\Phi_{h_1}^{-1}(g_2,h_2))
$$
identifies it with $X_G\times G$. Its output map is $(h_1,g_3,h_3,h)\mapsto(g_3,h_3)$, so the relative class is $[G][\varphi_G]$.

For part (b), write $P=X^0_{H_1,H_2}\times G$. The quotient $P\to\widetilde X^0_{H_1,H_2}$ is an $H_1$-torsor. Therefore, by Lemma \ref{lma:pullback_commutes_with_quot}, the fiber product used to compute the action of the operator is isomorphic to
$$
\widetilde X^0_{H_1,H_2}
\times_{\comvar}R_G(\Pi_{n,m})
\cong
\bigl(P\times_{\comvar}R_G(\Pi_{n,m})\bigr)/H_1.
$$

Before quotienting, the pullback consists of tuples $(h_1,g_2,h_2,g,g_4,h_4)$ satisfying
$$
(h_1,g_2,h_2)\in X^0_{H_1,H_2},\quad
[g_4,h_4]=e,\quad
g_4^nh_4^m=gh_1^mg^{-1}.
$$
The induced $H_1$-action is
$$
h\cdot(x,g,g_4,h_4)
=(h\cdot x,gh^{-1},g_4,h_4),
$$
where $x=(h_1,g_2,h_2)$ and $h\cdot x$ denotes simultaneous conjugation.

Consider the isomorphism
\begin{align*}
\Theta:P\times_{\comvar}R_G(\Pi_{n,m})
&\longrightarrow P\times H_2,\\
(h_1,g_2,h_2,g,g_4,h_4)
&\longmapsto
(h_1,g^{-1}g_4g,g^{-1}h_4g,g,
g_2^bh_2^{-a}),
\end{align*}
where $an+bm=1$.
The new commuting pair belongs to $F(h_1)$, and the last coordinate belongs to $H_2$. The inverse sends $(h_1,u,v,g,h')$ to
$$
\bigl(h_1,(h')^mh_1^{am},
(h')^{-n}h_1^{bm},g,gug^{-1},gvg^{-1}\bigr).
$$

The map $\Theta$ is $H_1$-equivariant for the actions described above, including conjugation on the extra $H_2$-factor. It preserves the output map, which on $P\times H_2$ ignores the last coordinate and agrees with $\varphi^0_{H_1,H_2}$ on $P$. Passing to the geometric quotients proves \eqref{eq:associated-centralizer-operator}.
\end{proof}

\begin{coro}\label{coro:associated-centralizer-eigenvector}
Under the hypotheses of Theorem \ref{tma:resultados}(b), suppose that
$$
[\pi:\mathcal E_{H_1,H_2}
\to\widetilde X^0_{H_1,H_2}]
=[H_2]\Unit{\widetilde X^0_{H_1,H_2}}
$$
in $\K\Varrel{\widetilde X^0_{H_1,H_2}}$.
Then
$$
\mathcal Z(\Pi_{n,m})
([\widetilde\varphi^0_{H_1,H_2}])
=[H_2][\widetilde\varphi^0_{H_1,H_2}].
$$
\end{coro}

\begin{proof}
By Theorem \ref{tma:resultados}(b) and the linearity of the pushforward $(\widetilde\varphi^0_{H_1,H_2})_!$, we obtain
\begin{align*}
\mathcal Z(\Pi_{n,m})
([\widetilde\varphi^0_{H_1,H_2}])
&=(\widetilde\varphi^0_{H_1,H_2})_!([\pi])\\
&=[H_2](\widetilde\varphi^0_{H_1,H_2})_!
(\Unit{\widetilde X^0_{H_1,H_2}})\\
&=[H_2][\widetilde\varphi^0_{H_1,H_2}],
\end{align*}
as claimed.
\end{proof}

\begin{obs}
The hypothesis of Corollary \ref{coro:associated-centralizer-eigenvector} holds, in particular, if $H_1$ is special or if $H_1$ centralizes $H_2$.

Indeed, if $H_1$ is special, the $H_1$-torsor $P\to\widetilde X^0_{H_1,H_2}$ is Zariski locally trivial. Hence, the associated bundle $\mathcal E_{H_1,H_2}$ is also Zariski locally trivial with fiber $H_2$, giving the required relative virtual class.

Alternatively, if $H_1$ centralizes $H_2$, its conjugation action on $H_2$ is trivial. Consequently,
$$
\mathcal E_{H_1,H_2}
\cong\widetilde X^0_{H_1,H_2}\times H_2,
$$
and the hypothesis follows directly.
\end{obs}

We now extend the preceding operator formula from the auxiliary varieties $\widetilde X^0_{H_1,H_2}$ to the original strata $X_{H_1,H_2}$. Let $$\widehat{\mathcal E}_{H_1,H_2}=(X^0_{H_1,H_2}\times G\times H_2)/M_{H_1,H_2}
\xrightarrow{\ \widehat\pi\ }\widehat X^0_{H_1,H_2},$$ where $r\in M_{H_1,H_2}$ conjugates $h_1,g_2,h_2,h'$ and sends $g$ to $gr^{-1}$.

\begin{coro}\label{coro:Nspecial2}
With the above notation,
\begin{equation}\label{eq:associated-bundle-operator}
\mathcal Z(\Pi_{n,m})([\varphi_{H_1,H_2}])
=[\widehat\varphi^0_{H_1,H_2}\circ\widehat\pi:
\widehat{\mathcal E}_{H_1,H_2}\to\comvar]_{\comvar}.
\end{equation}
In particular,
$$
\mathcal Z(\Pi_{n,m})([\varphi_{H_1,H_2}])
=[H_2][\varphi_{H_1,H_2}]
$$
whenever $\widehat\pi$ has relative class
$[H_2]\Unit{\widehat X^0_{H_1,H_2}}$.
\end{coro}

\begin{proof}
By Theorem \ref{tma:covering}, $[\varphi_{H_1,H_2}]=[\widehat\varphi^0_{H_1,H_2}]$. Apply Lemma \ref{lma:pullback_commutes_with_quot} to the $M_{H_1,H_2}$-torsor $X^0_{H_1,H_2}\times G\to\widehat X^0_{H_1,H_2}$. The isomorphism $\Theta$ in the preceding proof is also $M_{H_1,H_2}$-equivariant when $M_{H_1,H_2}$ conjugates the first three coordinates and the extra $H_2$-factor, obtaining \eqref{eq:associated-bundle-operator}.

The scalar formula follows by pushing forward the assumed relative class identity along $\widehat\varphi^0_{H_1,H_2}$.
\end{proof}

\begin{obs}
The hypothesis of Corollary \ref{coro:Nspecial2} holds, in particular, under either of the following conditions:
\begin{itemize}
    \item $M_{H_1,H_2}$ is special;
    \item $N_{H_1,H_2}$ is special and either
    $H_1$ is special or $H_1$ centralizes $H_2$.
\end{itemize}

To see this, consider the following Cartesian diagram:
\[
\begin{tikzcd}[column sep=large, row sep=large]
\mathcal E_{H_1,H_2}
    \arrow[r]
    \arrow[d, "\pi"']
&
\widehat{\mathcal E}_{H_1,H_2}
    \arrow[d, "\widehat\pi"]
\\
\widetilde X^0_{H_1,H_2}
    \arrow[r, "q"']
&
\widehat X^0_{H_1,H_2},
\end{tikzcd}
\]
where $q$ is the $N_{H_1,H_2}$-torsor of Proposition \ref{prop:quotient-torsor}. Thus, $\mathcal E_{H_1,H_2}$ is the pullback of $\widehat{\mathcal E}_{H_1,H_2}$ along $q$.

If $M_{H_1,H_2}$ is special, the associated bundle $\widehat\pi$ is Zariski locally trivial with fiber $H_2$. For the second condition, if $H_1$ is special, the quotient map $X^0_{H_1,H_2}\times G\times H_2 \to\mathcal E_{H_1,H_2}$ is a Zariski locally trivial $H_1$-torsor. In particular, the $H_1$-torsor $X^0_{H_1,H_2}\times G\to\widetilde X^0_{H_1,H_2}$ is also Zariski locally trivial, and its local sections trivialize the associated bundle $\pi$ with fiber $H_2$. Alternatively, if $H_1$ centralizes $H_2$, the bundle $\pi$ is a product bundle. Since $N_{H_1,H_2}$ is special, the torsor $q$ admits Zariski local sections. By the Cartesian diagram, these sections allow the local trivializations of $\pi$ to descend to local trivializations of $\widehat\pi$.

Thus, in either case,
$$
[\widehat\pi]
=[H_2]\Unit{\widehat X^0_{H_1,H_2}},
$$
as required.
\end{obs}

\begin{obs}\label{rem:general-algebraic-groups}
Although we have assumed that $G$ is reductive, the results of this section extend to arbitrary algebraic groups, with some minor technical adjustments. Recall that the notation $/$ refers to geometric quotients of free group actions, rather than affine GIT quotients. In our setting, their existence follows from the explicit constructions above, without requiring reductivity.

First, the loci $X^0_{H_1,H_2}$ and $X_{H_1,H_2}$ are constructible for any pair of closed algebraic subgroups $H_1\leq H_2\leq G$. Indeed, for any closed subgroup $H\leq G$, the condition $\Stab_G(h)=H$ is constructible. The inclusion $H\subseteq \Stab_G(h)$ defines a closed condition, while the failure of the reverse inclusion is equivalent to the existence of an element in $G\setminus H$ commuting with $h$. By Chevalley's theorem \cite[Theorem 29.23.3, Tag 054K]{stacks-project}, the latter condition is constructible. Applying this to $h_1$ and $h_1^m$ shows that $X^0_{H_1,H_2}$ is constructible. Moreover, $X_{H_1,H_2}$ is the image of $X^0_{H_1,H_2}\times G$ under simultaneous conjugation, and is therefore constructible by Chevalley's theorem. However, finiteness of the stratification \eqref{eq:stratification} is no longer guaranteed and must be imposed whenever the decomposition into a finite sum of relative virtual classes is used.

Second, since the subgroups $H_1$, $H_2$ and $M_{H_1,H_2}=N_G(H_1)\cap N_G(H_2)$ are closed, the corresponding homogeneous spaces exist \cite[Theorem B.37]{Milne2017}. Consequently, the geometric quotients $\widetilde X^0_{H_1,H_2}$ and $\widehat X^0_{H_1,H_2}$ can be constructed as above, without requiring reductivity.

Thus, provided that the stratification has finitely many simultaneous conjugacy classes of stabilizer pairs, the quotient constructions and operator formulas established in this section remain valid for arbitrary algebraic groups $G$. Moreover, we illustrate the method with explicit computations for the non-reductive complex algebraic group $\AGL_1(\mathbb C)$ in Section \ref{sec:AGL1} and the reductive complex algebraic group $\SL_2(\mathbb C)$ in Section \ref{sec:SL2}.
\end{obs}

\section{\texorpdfstring{$\Aff{1}{\mathbb{C}}$}{AGL(1,C)}-Representation Varieties of Torus Links}\label{sec:AGL1}

In this section we show explicit computations of the virtual classes of representation varieties of torus links with the quantum method developed in Section \ref{sec:Quantum_Method}. We consider the group of affine transformations of the affine line over $\C$,
$$\AGL_1(\C) = \left\{ \begin{pmatrix} a & b \\ 0 & 1 \end{pmatrix} : a \in \C^*, b \in \C \right\},$$
where the group operation is given by matrix multiplication. We abbreviate the notation by $\AGL_1$. It easily follows from the definition that the virtual class of $\AGL_1$ is $[\AGL_1]=\lL(\lL-1)$. The $\AGL_1$-representation variety of an $(n,m)$-torus link with $d$ components is given by
\begin{equation*}
R_G(L^d_{n,m})=
\left\{
(A,B,F_1,\ldots,F_{d-1})\in\AGL_1(\C)^{d+1}
\ \middle|\
\begin{aligned}
A^n &= B^m,\\
F_iA^nF_i^{-1} &= A^n,\quad i=1,\ldots,d-1
\end{aligned}
\right\}.
\end{equation*}

In light of the results of Section \ref{sec:QuantizationVirtualClasses}, we first compute the operator $\mathcal{Z}(\Pi_{n,m}) : \K\Varrel{\comvar} \to \K\Varrel{\comvar} $. Applying topological recursion, it suffices to calculate $\mathcal{Z}(\Pi_{n,m})^{d}$ and compose with the operators $\mathcal{Z}(\mathbf{T})$ and $\mathcal{Z}(\mathbf{T}^\dagger)$ to obtain the virtual class $[R_{G}(L^d_{n,m})]$. 

As done in Section \ref{sec:computation_virtualclass}, we shall stratify the variety $X=G\times_{\comvar}R_G(\Pi_{n,m})=\{(h_1,g_2,h_2)\in G^3 \ | \ [g_2,h_2]=e, h_1^m=g_2^{n}h_2^{m}\}$ considering the stratification given in Equation (\ref{eq:stratification}) in order to find the eigenvectors of $\mathcal{Z}(\Pi_{n,m})$. To this end, we study the stabilizers of $h_1$ and $h_1^m$ in $\AGL_1$. Let us write
\[
U = \left\{ \begin{pmatrix} 1 & y \\ 0 & 1 \end{pmatrix} \in G \mathrel{\bigg|} y \in \mathbb{C} \right\} \cong \C, \qquad T = \left\{ \begin{pmatrix} x & 0 \\ 0 & 1 \end{pmatrix} \in G \mathrel{\bigg|} x \in \mathbb{C}^* \right\} \cong \C^*.
\]

These subgroups are the stabilizers under conjugation of nonidentity unipotent and nonidentity diagonal matrices, respectively. Together with $G$, which is the stabilizer of the identity, they account for all possible stabilizers up to conjugation, since every matrix with $a\neq 1$ is conjugate to an element of $T$. 

Without loss of generality, we may assume that $h_1$ belongs to $U$ or $T$, up to conjugation. Indeed, if its upper-left entry is $a\neq 1$, then $h_1$ is diagonalizable and therefore conjugate to an element of $T$.

\begin{itemize}
    \item $(H_1,H_2)=(\AGL_1,\AGL_1)$. The only case in which the stabilizer is the whole group occurs when $h_1=h_1^m=I$.
    \item $(H_1,H_2)=(U,U)$. This case corresponds to $h_1^m$ being a unipotent matrix with $y\neq 0$. Consequently, $h_1$ is a unipotent matrix.
    \item $H_1=T$. Considering the stabilizer of $h_1$ as the maximal torus corresponds to $h_1\in T$ with $x\neq 1$. We shall consider different cases:
        \begin{itemize}
            \item $h_1^m = I$. This is the case when $x$ is an $m$-th root of unity. The corresponding pair is $(H_1,H_2)=(T,\AGL_1)$.
            \item $h_1^m \neq  I$. In this case $h_1^m$ is also in $T$ with $x^m \neq 1$, hence $(H_1,H_2)=(T,T)$.
        \end{itemize}
\end{itemize}

This provides a complete characterization of the stabilizer pairs, yielding the stratification
\begin{equation}
\label{eq:stratification_AGL1}
X=X_{\AGL_1,\AGL_1}\sqcup X_{U,U}
\sqcup X_{T,\AGL_1}\sqcup X_{T,T}.
\end{equation}
Our goal is to identify a submodule invariant under $\mathcal{Z}(\Pi_{n,m})$, allowing the operator to be expressed in matrix form so the topological recursion argument can be applied.

The first eigenvector is a direct consequence of Theorem \ref{tma:resultados}(a). Considering $H_2=\AGL_1$, it directly follows that $[\varphi_G:X_G \to \comvar]$ with $X_G:=X_{G,G} \sqcup X_{T,G}$ is an eigenvector with eigenvalue $\lL(\lL-1)$.

For the pair $(T,T)$, $N_G(T)/T$ is trivial and $T$ is abelian, so $H_1$ centralizes $H_2$. Therefore, we apply Corollary \ref{coro:Nspecial2} to obtain that $[\varphi_{T,T}:X_{T,T} \to \comvar]$ is an eigenvector with eigenvalue $[T]=[\C^*]=\lL-1$. Finally, for the pair $(U,U)$ we have $N_{G}(U)/U=\C^*$, which is special, and $U$ is abelian. Thus, again by Corollary \ref{coro:Nspecial2}, the class $[\varphi_{U,U}:X_{U,U} \to \comvar]$ is an eigenvector of $\mathcal{Z}(\Pi_{n,m})$ with eigenvalue $[U]=\lL$.

For $G=\AGL_1(\C)$, we summarize the obtained results in the following proposition:

\begin{prop}The images of the virtual classes defined above under the operator $\mathcal{Z}(\Pi_{n,m}) : \K\Varrel{\comvar} \to \K\Varrel{\comvar}$ are the following:
\begin{align*}
    \mathcal{Z}(\Pi_{n,m})([\varphi_G:X_G \to \comvar]) &= \lL(\lL-1)[\varphi_G:X_G \to \comvar] \\
    \mathcal{Z}(\Pi_{n,m})([\varphi_{T,T}:X_{T,T} \to \comvar]) &= (\lL-1)[\varphi_{T,T}:X_{T,T} \to \comvar] \\
    \mathcal{Z}(\Pi_{n,m})([\varphi_{U,U}:X_{U,U} \to \comvar]) &= \lL[\varphi_{U,U}:X_{U,U} \to \comvar]
\end{align*}
Thus, the submodule of $\K\Varrel{\comvar}$ spanned by them is invariant under $\mathcal{Z}(\Pi_{n,m})$.
\end{prop}

\begin{obs}\label{obs:matrixtoriAGL1}
Considering the invariant submodule spanned by the virtual classes, the operator $\mathcal{Z}(\Pi_{n,m})$ can be represented, in the basis $\{[\varphi_G],[\varphi_{T,T}],[\varphi_{U,U}]\}$ described above, by the matrix
\begin{equation}\label{eq:matrix_agl1}
\mathcal{Z}(\Pi_{n,m})= \left(\begin{matrix}
\lL(\lL-1) & 0 & 0  \\
0 & \lL-1 & 0  \\
0 & 0 & \lL  
\end{matrix}\right)
\end{equation}
\end{obs}

Since $X=X_G\sqcup X_{T,T}\sqcup X_{U,U}$, we have
$$
[\varphi]=[\varphi_G]+[\varphi_{T,T}]+[\varphi_{U,U}]
$$
in $\K\Varrel{\comvar}$. By Theorem \ref{tma:TQFTcomptesVirtClass}, the virtual class of the representation variety is computed as
$$
[R_G(L^d_{n,m})]
=\mathcal{Z}(E(\sS^3,L_{n,m}^d))(\Unit{\star})
=\mathcal{Z}(\mathbf{T}^\dagger)\circ
\mathcal{Z}(\Pi_{n,m})^d\circ
\mathcal{Z}(\mathbf{T})(\Unit{\star}).
$$
Since $[\varphi:X\to\comvar] =\mathcal{Z}(\Pi_{n,m})\circ \mathcal{Z}(\mathbf{T})(\Unit{\star})$, the remaining iteration on these virtual classes has exponent $d-1$. Accordingly, we shall compute the images of the aforementioned virtual classes under the operator $\mathcal{Z}(\mathbf{T}^\dagger): \K\Varrel{\comvar}\to\K{\Var{\C}}$. Therefore, we need to compute the virtual classes of the corresponding strata in $R_G(L^1_{n,m})$.

\begin{prop}\label{prop:cerrar_agl1}
The images of the virtual classes under the operator $\mathcal{Z}(\mathbf{T}^\dagger)$ for $G=\AGL_1(\C)$ are the following:
\begin{align*}
    \mathcal{Z}(\mathbf{T}^\dagger)([\varphi_G:X_G \to \comvar]) &=\left(1+\left(m-1\right)\lL\right)\left(1+\left(n-1\right)\lL\right) \\
    \mathcal{Z}(\mathbf{T}^\dagger)([\varphi_{T,T}:X_{T,T} \to \comvar]) &= \lL\left(\lL-1-nm\right) \\
    \mathcal{Z}(\mathbf{T}^\dagger)([\varphi_{U,U}:X_{U,U} \to \comvar]) &= \lL-1
\end{align*}
\end{prop}

\begin{proof}
Applying the operator involves the composition of the spans for each stratum $X=X_G \sqcup X_{T,T} \sqcup X_{U,U}$. We explicitly compute the composition for $X_G$. The remaining strata are treated in an analogous way. The pullback is described by the following diagram:

\[
\begin{tikzcd}[
    column sep={2cm, between origins}, 
    row sep={1.5cm, between origins}, 
    ampersand replacement=\&
]
\& \& X_G \times_{\comvar} G \arrow[dl, "t_1"'] \arrow[dr, "t_2"] \& \& \\
\& X_G \arrow[dl, "c"'] \arrow[dr, "\varphi_G"] \& \& G \arrow[dl, "q_1"'] \arrow[dr, "d"] \& \\
\star \& \& \comvar \& \& \star
\end{tikzcd}
\]
where $\varphi_G:X_G \to \comvar$ is given by the restriction of $(h_1,g_2,h_2)\mapsto (g_2,h_2)$, and $q_1:G\to\comvar$ is given by $g\mapsto(g,e)$. The resulting variety can be described as $Y_G:=X_G\times_{\comvar}G=\{(h,\bar{h}) \in G^2 : h^m=\bar{h}^n, \Stab(h^m)_G=G \}$. The same pullback for the remaining strata gives
\begin{align*}
Y_{T,T}&:=\{(h,\bar{h})\in G^2 \mid
h^m=\bar{h}^n,\ h^m\in G\setminus U\},\\
Y_{U,U}&:=\{(h,\bar{h})\in G^2 \mid
h^m=\bar{h}^n,\ h^m\in U\setminus\{I\}\}.
\end{align*}
These conditions correspond to the stabilizer pairs $(T,T)$ and $(U,U)$, respectively, up to conjugation. The rest of the proof consists of computing the virtual classes of $Y_G$, $Y_{T,T}$ and $Y_{U,U}$.

For the stratum $Y_G$, the condition $\Stab_G(h^m)=G$ implies $h^m=\bar{h}^n=I$. We identify each element $h\in \AGL_1$ with the pair $(x,y)\in\mathbb C^*\times\mathbb C$. The requirement $h^m=(x^m , y\sum_{k=0}^{m-1}x^k)=(1,0)$ implies $x^m=1$ and $y\sum_{k=0}^{m-1}x^k=0$. Hence, if $x=1$ then $y=0$, yielding the identity element. If $x\in \mu_m\setminus\{1\}$, then $y \in \mathbb{C}$ for each of the $m-1$ nontrivial roots of unity. The same reasoning applies to $\bar{h}^n=I$. Since $Y_G=\{ h\in G \ | \ h^m=I \} \times \{ \bar{h} \in G \ | \ \bar{h}^n=I \}$, the virtual class is $[Y_G]=(1+(m-1)\lL)(1+(n-1)\lL)$.

For the unipotent stratum $Y_{U,U}$, the condition $\Stab_G(h^m)=U$ requires $h^m=\bar{h}^n \in U \setminus \{I\}$. In this case, $h^m = (1, my)$, hence there is exactly one pair $(h,\bar{h})$ for each nontrivial unipotent matrix, so $[Y_{U,U}]=\lL-1$.

We now consider $Y_{T,T}$, where the stabilizer is the maximal torus $T$. Writing $h=(x,y)$ and $\bar h=(\bar x,\bar y)$, the condition $h^m=\bar h^n$ implies that the multiplicative coordinates belong to the curve
$$
C:=\{(x,\bar x)\in\mathbb C^*\times\mathbb C^*
\mid x^m=\bar x^n\neq 1\}.
$$
Since $n$ and $m$ are coprime, the parametrization $t\mapsto(t^n,t^m)$ identifies $C$ with $\mathbb C^*$ minus the $nm$-th roots of unity. Hence, $[C]=\lL-1-nm$.

The additive components are determined by one parameter in the affine line. More explicitly, the equality $h^m=\bar h^n$ imposes the additional equation
$$
y\sum_{k=0}^{m-1}x^k
=\bar y\sum_{k=0}^{n-1}\bar x^k.
$$
Since $x^m=\bar x^n\neq 1$, both coefficients are nonzero. Thus, for each fixed pair $(x,\bar x)\in C$, the set of pairs $(y,\bar y)$ is isomorphic to $\mathbb A^1$. Finally, the virtual class is given by
$$[Y_{T,T}]=\lL[C]=\lL(\lL-1-nm).$$
\end{proof}

Thanks to Theorem \ref{tma:TQFTcomptesVirtClass} we can use the aforementioned stratification and Proposition \ref{prop:cerrar_agl1} to compute the virtual class of the representation variety with the matrix product $\mathbf{w}\cdot \mathcal{Z}(\Pi_{n,m})^{d-1}\cdot \mathbf{v}$, where $\mathbf{v}$ is the column vector whose three entries are equal to $1$ and $\mathbf{w}$ is given by the components calculated in Proposition \ref{prop:cerrar_agl1}. The matrix of the operator $\mathcal{Z}(\Pi_{n,m})$ is already diagonal, hence
\begin{equation*}
\mathbf{w} \cdot \mathcal{Z}(\Pi_{n,m})^{d-1} \cdot \mathbf{v}= \begin{pmatrix}
(1+(m-1)\lL)(1+(n-1)\lL) \\
\lL(\lL-1-mn) \\
\lL-1
\end{pmatrix}^{T} \left(\begin{matrix}
\lL(\lL-1) & 0 & 0  \\
0 & \lL-1 & 0  \\
0 & 0 & \lL  
\end{matrix}\right)^{d-1} \begin{pmatrix}
    1 \\ 1 \\ 1
\end{pmatrix}
\end{equation*}

Performing the calculation, we have obtained the following result:

\begin{tma}\label{tma:virtual_class_AGL1}
For $G=\AGL_1(\C)$ and $d,n,m\geq 1$, with $n$ and $m$ coprime, the virtual class of the torus link representation variety is
\begin{align*}
[R_{\AGL_1(\mathbb{C})}(L_{n,m}^d)] &= \lL(\lL-1-nm)(\lL-1)^{d-1} \\
&\quad + (1+(m-1)\lL)(1+(n-1)\lL)(\lL(\lL-1))^{d-1} \\
&\quad + (\lL-1)\lL^{d-1}.
\end{align*}
\end{tma}

As a corollary, we recover the virtual class of the $\AGL_1(\C)$-representation variety of an $(n,m)$-torus knot, previously computed in \cite{gplogmunlowrankAGL1}.

\begin{coro}
    The virtual class of the $\AGL_1(\C)$-representation variety of an $(n,m)$-torus knot is 
    \begin{equation*}
         [R_{\AGL_1(\mathbb{C})}(K_{n,m})]=(mn-m-n + 2)(\lL^2-\lL).
    \end{equation*}
\end{coro}

\begin{proof}
    The result follows by evaluating the virtual class in Theorem \ref{tma:virtual_class_AGL1} at $d=1$. In this case, the product $\mathbf{w}\cdot \mathcal{Z}(\Pi_{n,m})^{0}\cdot \mathbf{v}$ amounts to summing the contributions of the strata, as described in Proposition \ref{prop:cerrar_agl1}.
\end{proof}

\section{\texorpdfstring{$\SL_2(\mathbb{C})$}{SL(2,C)}-Representation Varieties of Torus Links}\label{sec:SL2}

In this section we apply the quantum method to compute the virtual class of the representation variety of an $(n,m)$-torus link with $d$ components for the group $\SL_2(\C)$, which we abbreviate by $\SL_2$. We proceed exactly as in the previous section to determine the operator $\mathcal{Z}(\Pi_{n,m})$. Again, we shall stratify the variety $X=G\times_{\comvar}R_G(\Pi_{n,m})=\{(h_1,g_2,h_2)\in G^3 \ | \ [g_2,h_2]=e, h_1^m=g_2^{n}h_2^{m}\}$ considering the stratification given in Equation (\ref{eq:stratification}).

Let us write
\[
U = \left\{ \pm \begin{pmatrix} 1 & y \\ 0 & 1 \end{pmatrix} \in G \mathrel{\bigg|} y \in \mathbb{C} \right\}, \qquad T = \left\{ \begin{pmatrix} x & 0 \\ 0 & x^{-1} \end{pmatrix} \in G \mathrel{\bigg|} x \in \mathbb{C}^* \right\}.
\]

These subgroups are the stabilizers under conjugation of noncentral elements of trace $\pm 2$ and regular semisimple elements, respectively. Note that, unlike in Section~\ref{sec:AGL1}, the group $U$ is not connected: it is the direct product of the center $\{\pm I\}$ of $\SL_2$ and the unipotent subgroup, so that $[U]=2\lL$. Together with $G$, which is the stabilizer of the central elements $\pm I$, they account for all possible stabilizers up to conjugation.

\begin{itemize}
    \item $(H_1,H_2)=(\SL_2,\SL_2)$.
    The stabilizer is the whole group precisely when $h_1=\pm I$, and hence $h_1^m=\pm I$.

    \item $(H_1,H_2)=(U,U)$. This case corresponds to $h_1^m$ being a noncentral matrix of trace $\pm 2$, that is, a nontrivial unipotent element or its negative. Consequently, $h_1$ is also of this type.

    \item $H_1=T$. In this case, $h_1$ is regular semisimple, with eigenvalues $\lambda$ and $\lambda^{-1}$, where $\lambda\neq\pm 1$. We distinguish two cases:
    \begin{itemize}
        \item $h_1^m=\pm I$. This occurs when $\lambda^{2m}=1$, giving $(H_1,H_2)=(T,\SL_2)$.
        \item $h_1^m\neq\pm I$. In this case, $h_1^m$ is also regular semisimple, with the same centralizer $T$. Hence $(H_1,H_2)=(T,T)$.
    \end{itemize}
\end{itemize}

This provides a complete characterization of the stabilizer pairs, yielding the stratification
\begin{equation}\label{eq:stratification_SL2}
X=X_{\SL_2,\SL_2}\sqcup X_{U,U} \sqcup X_{T,\SL_2}\sqcup X_{T,T}.
\end{equation}
Equivalently, setting $X_G:=X_{\SL_2,\SL_2}\sqcup X_{T,\SL_2}$, we have
$$ X=X_G\sqcup X_{T,T}\sqcup X_{U,U}.$$

Before proceeding with the analysis of the strata, we focus our attention on the stratum $(T,T)$. Indeed, the Weyl involution, that is, the $\Z_2$-action $t\mapsto t^{-1}$, acts by permuting the eigenvalues of the corresponding matrix elements, so the ordinary virtual class does not retain all the information required for this stratum. We therefore work with equivariant virtual classes, mentioned in Section \ref{sec:groth_ring}. We collect in the following lemma some virtual classes and their equivariant refinements that will be useful in subsequent computations.

\begin{lema}\label{lem:some_epoly}
The virtual classes of the following subvarieties of $\SL_2(\mathbb{C})$ are given by:

\begin{enumerate}[label=(\roman*)]
\item The conjugacy classes $\mathcal O_{J_\pm}$, where $J_+$ is a nontrivial unipotent element and $J_-=-J_+$, of traces $2$ and $-2$, respectively. Both have virtual class $[\mathcal O_{J_\pm}]=\lL^2-1$.
\item The regular semisimple orbit: $[\SL_2/T]=\lL^2+\lL$.
\item For the Weyl involution on the maximal torus $T$, the coefficients of the equivariant class are $[T]^+=\lL$ and $[T]^-=-1$.
\item For the Weyl involution on $\SL_2/T$, the coefficients of the equivariant class are $[\SL_2/T]^+=\lL^2$ and $[\SL_2/T]^-=\lL$.
\end{enumerate}
\end{lema}

\begin{proof}
We refer the reader to \cite[Proposition 3.2(i),(iii)]{lomune2013surfaces} for a detailed derivation of the first and second statements. The third assertion follows directly considering $T\cong\C^*$, so that the Weyl involution becomes $x\mapsto x^{-1}$. Therefore $[T]^{+}=[\C^* / \Z_2]=\lL$ and $[T]^{-}=[\C^*]-\lL=(\lL-1)-\lL=-1$.

For the last assertion, as in \cite[Proposition 3.2(ii)]{lomune2013surfaces}, identify $\SL_2/T$ with the variety of ordered pairs of distinct points of $\mathbb P^1$. The Weyl involution exchanges the two points, so its quotient is $\operatorname{Sym}^2(\mathbb P^1)\setminus\Delta\cong\mathbb P^2\setminus C$, where $C$ is a conic. Thus, $[\SL_2/T]^+=[\mathbb P^2]-[C]=\lL^2$, and $[\SL_2/T]^-=(\lL^2+\lL)-\lL^2=\lL$.
\end{proof}

We now resume the analysis of the stratification. Our aim is to determine suitable eigenvectors and eigenvalues for the operator $\mathcal{Z}(\Pi_{n,m})$. The first eigenvector follows directly from Theorem \ref{tma:resultados}(a). Taking $H_2=\SL_2$, it follows that $[\varphi_G:X_G \to \comvar]$ with $X_G:=X_{G,G} \sqcup X_{T,G}$ is an eigenvector with eigenvalue $[\SL_2]=\lL^3-\lL$. The computation of this virtual class follows as in \cite[Lemma 3.1]{lomune2013surfaces}. 

Secondly, the stratum $X_{U,U}$ gives an eigenvector by virtue of Corollary \ref{coro:Nspecial2}. Specifically, $N_{U,U}=N_{\SL_2}(U)/U\cong \C^*$, which is a connected reductive special group, and $U$ is abelian, so $H_1$ centralizes $H_2$. Hence, $[\varphi_{U,U}:X_{U,U} \to \comvar]$ is an eigenvector with eigenvalue $[U]=2\lL$.

The analysis of the final stratum, $X_{T,T}$, requires a more careful treatment. As in Section \ref{sec:QuantizationVirtualClasses}, consider the subvariety
\[
X_{T,T}^0:=\{(h_1,g_2,h_2)\in X \ | \ (\Stab_{\SL_2}(h_1),\Stab_{\SL_2}(h_1^m))= (T,T)\}.
\]

Let $\widetilde{X}_{T,T}^0:=(X^0_{T,T}\times \SL_2) / T$. Considering the morphism $\widetilde{\varphi}_{T,T}^0:\widetilde{X}_{T,T}^0\to \comvar$ and applying Corollary \ref{coro:associated-centralizer-eigenvector}, whose centralization hypothesis holds because $T$ is abelian, we obtain that its virtual class is an eigenvector of $\mathcal{Z}(\Pi_{n,m})$, with eigenvalue $[T]=\lL-1$. However, unlike the unipotent case, we cannot simply invoke Corollary \ref{coro:Nspecial2} to claim that $[\varphi_{T,T}:X_{T,T}\to \comvar]$ is an eigenvector. This phenomenon is governed by the Weyl involution arising from the normalizer quotient $N_{T,T}=N_{\SL_2}(T)/T \cong \Z_2$. To analyze this case, we use the equivariant virtual classes introduced previously.

Let us denote by $[\widetilde{\varphi}^0_{T,T}]^\pm$ the $\mathbb Z_2$-equivariant components of $[\widetilde{\varphi}_{T,T}^0:\widetilde{X}_{T,T}^0\to \comvar]$ with respect to the action of $N_{T,T} \cong \mathbb Z_2$. Since $T$ acts trivially on $X^0_{T,T}$, as all its elements are diagonal, we have 
\begin{equation}\label{eq:XT0_prod}
\widetilde X_{T,T}^0=(X_{T,T}^0\times \SL_2)/T \cong X_{T,T}^0\times (\SL_2/T).
\end{equation}

The Weyl group acts diagonally on this product, and its action is free on the factor $\SL_2/T$. By Theorem \ref{tma:covering} the quotient $\widetilde X_{T,T}^0 / N_{T,T}$ is identified, at the level of relative virtual classes, with $X_{T,T}$. Hence, $[\varphi_{T,T}]=[\widetilde{\varphi}^0_{T,T}]^+$. Although Corollary \ref{coro:Nspecial2} does not establish that $[\varphi_{T,T}]$ is an eigenvector, Equation (\ref{eq:associated-bundle-operator}) still provides an expression for its image under $\mathcal Z(\Pi_{n,m})$. Indeed, the associated bundle is $\widehat{\mathcal{E}}_{T,T}=(X^0_{T,T}\times \SL_2 \times T)/ N_G(T) \cong (T\times \widetilde{X}^0_{T,T})/\mathbb Z_2$. Hence, we have
$$\mathcal Z(\Pi_{n,m})\bigl([\varphi_{T,T}]\bigr)=[T\times\widetilde X^{0}_{T,T}]^{+}.$$

The multiplicativity hypothesis of Lemma \ref{lem:plusminusformula} holds, since the class of $T$ with the Weyl involution belongs to $V_{\mathbb{Z}_2}$ (see \cite[Example~3.6.22]{vogelthesis} for further details), and Lemma \ref{lem:plusminusformula} together with Lemma \ref{lem:some_epoly}(iii) gives
$$[T\times\widetilde X^{0}_{T,T}]^{+}=[T]^{+}[\widetilde{\varphi}^0_{T,T}]^{+}+[T]^{-}[\widetilde{\varphi}^0_{T,T}]^{-}=(\lL+1)\,[\widetilde{\varphi}^0_{T,T}]^{+}-[\widetilde{\varphi}^0_{T,T}].$$

Putting together all the computations, we have proven the following result:

\begin{prop}The images of the virtual classes defined above under the operator $\mathcal{Z}(\Pi_{n,m}) : \K\Varrel{\comvar} \to \K\Varrel{\comvar}$ are the following:
\begin{align*}
    \mathcal{Z}(\Pi_{n,m})([\varphi_G:X_G \to \comvar]) &= (\lL^3-\lL)[\varphi_G:X_G \to \comvar] \\
    \mathcal{Z}(\Pi_{n,m})([\varphi_{T,T}:X_{T,T} \to \comvar]) &= (\lL+1)[\varphi_{T,T}:X_{T,T} \to \comvar] - [\widetilde{\varphi}_{T,T}^0:\widetilde{X}_{T,T}^0\to \comvar]  \\
    \mathcal{Z}(\Pi_{n,m})([\widetilde{\varphi}_{T,T}^0:\widetilde{X}_{T,T}^0\to \comvar]) &=(\lL-1)[\widetilde{\varphi}_{T,T}^0:\widetilde{X}_{T,T}^0\to \comvar] \\
    \mathcal{Z}(\Pi_{n,m})([\varphi_{U,U}:X_{U,U} \to \comvar]) &= 2\lL[\varphi_{U,U}:X_{U,U} \to \comvar]
\end{align*}

Thus, the submodule of $\K\Varrel{\comvar}$ spanned by them is invariant under $\mathcal{Z}(\Pi_{n,m})$.
\end{prop}

\begin{obs}\label{obs:matrixtoriSL2}
   In line with Remark \ref{obs:matrixtoriAGL1}, the operator $\mathcal{Z}(\Pi_{n,m})$ for $G=\SL_2$ is represented, in the basis $\{[\varphi_G],[\varphi_{T,T}], [\widetilde{\varphi}^0_{T,T}],[\varphi_{U,U}]\}$ described above, by the matrix

\begin{equation}\label{eq:matrix_sl2}
\mathcal{Z}(\Pi_{n,m})= \left(\begin{matrix}
\lL^3-\lL & 0 & 0 & 0  \\
0 & \lL+1 & 0  & 0\\
0 & -1 & \lL-1  & 0 \\
0 & 0 & 0 & 2\lL
\end{matrix}\right)
\end{equation}
\end{obs}

As in the $\AGL_1$ case, we have $X=X_G\sqcup X_{T,T}\sqcup X_{U,U}$, and hence
$$
[\varphi]=[\varphi_G]+[\varphi_{T,T}]+[\varphi_{U,U}].
$$
However, the additional class $[\widetilde{\varphi}^0_{T,T}]$ is required to obtain an invariant submodule under $\mathcal{Z}(\Pi_{n,m})$.

\begin{prop}\label{prop:cerrar_sl2}
Let $m$ be odd. The images of the virtual classes under the operator $\mathcal{Z}(\mathbf{T}^\dagger)$ for $G=\SL_2(\C)$ are the following:
\begin{align*}
    \mathcal{Z}(\mathbf{T}^\dagger)([\varphi_G:X_G \to \comvar]) &=\left(1 + \frac{m-1}{2}(\lL^2+\lL)\right) \left(2 + (n-1)(\lL^2+\lL)\right) \\
    \mathcal{Z}(\mathbf{T}^\dagger)([\varphi_{T,T}:X_{T,T} \to \comvar]) &= \lL^3-(mn+1)\lL^2-mn\lL \\
    \mathcal{Z}(\mathbf{T}^\dagger)([\widetilde{\varphi}_{T,T}^0:\widetilde{X}_{T,T}^0\to \comvar]) &= (\lL-1-2mn)(\lL^2+\lL) \\
    \mathcal{Z}(\mathbf{T}^\dagger)([\varphi_{U,U}:X_{U,U} \to \comvar]) &= 2(\lL^2-1)
\end{align*}

\end{prop}

\begin{proof}
Proceeding as in the proof of Proposition \ref{prop:cerrar_agl1}, we compute the image of the virtual classes under the operator $\mathcal{Z}(\mathbf{T}^\dagger)$. The fact that the stabilizers coincide with the entire group implies $h^m=\bar{h}^n=\pm I$. Consequently, we define $Y_G:=X_G \times_{\comvar} G$, which decomposes as the disjoint union: 
$$\{(h, \bar{h}) \mid h^m = \bar{h}^n = I\} \sqcup \{(h, \bar{h}) \mid h^m = \bar{h}^n = -I\}.$$
We define the set of $m$-th roots of an element $c\in G$ as $S_m(c)=\{h\in G \ | \ h^m=c \}$. The virtual class is computed as 
$$[Y_G]=[S_m(I)]\cdot[S_n(I)]+[S_m(-I)]\cdot[S_n(-I)].$$ 
Since $m$ is assumed to be odd, the solutions to the equations $h^m = \pm I$ consist of the corresponding central elements together with a collection of regular semisimple elements whose eigenvalues are distinct from $\pm 1$. The nontrivial $m$-th roots of unity form $(m-1)/2$ pairs under inversion, each indexing a regular semisimple conjugacy class in $\SL_2$. Considering changes of basis in $\SL_2$, each of these classes is isomorphic to $G/T$. Therefore,
    \[
    [S_m(I)]=[S_m(-I)]=1+\frac{m-1}{2}(\lL^2+\lL).
    \]

    For the other part of the disjoint union, $\bar h$ must satisfy $\bar h^n=\pm I$. We consider the combined solution set, which contains the two central elements $\pm I$ and $2n-2$ regular roots. Considering the $n-1$ conjugacy classes determined by the remaining eigenvalues, paired under inversion, we have
    \[
    [S_n(I)]+[S_n(-I)]=2+(n-1)(\lL^2+\lL).
    \]

    We remark that the result is independent of the parity of $n$. By taking the product of both results, we obtain the virtual class of $Y_G$:
    \[
    [Y_G] = \left(1 + \frac{m-1}{2}(\lL^2+\lL)\right) \left(2 + (n-1)(\lL^2+\lL)\right).
    \]

    The next stratum we consider is $Y_{U,U}$. In this case $h^m=\bar{h}^n=c$ is noncentral and unipotent up to sign. The set of noncentral matrices with trace $\pm 2$ is partitioned into two orbits, $\mathcal{O}_{J_+}$ and $\mathcal{O}_{J_-}$, corresponding to traces $2$ and $-2$, respectively.

    Since $m$ is assumed to be odd, the equation $h^m=c$ admits a unique solution unipotent up to sign having the same trace as $c$. For $\bar{h}^n$, the number of solutions depends on the parity of $n$.

    \begin{itemize}
        \item If $n$ is odd, following the argument above for $h$ gives a unique solution unipotent up to sign for $\bar{h}^n=c$.
        \item If $n$ is even, the equation yields two solutions if $c\in \mathcal{O}_{J_+}$ and no solutions if $c\in \mathcal{O}_{J_-}$, as the even power must have trace $2$.
    \end{itemize}

In view of the fact that $Y_{U,U}$ maps to $\mathcal{O}_{J_+}\sqcup \mathcal{O}_{J_-}$ via the map $(h,\bar{h}) \mapsto c$, we take the product of the number of solutions for $h$ and $\bar{h}$. If $n$ is odd, there is one solution pair $(h,\bar{h})$ over each point of either $\mathcal{O}_{J_+}$ or $\mathcal{O}_{J_-}$. According to Lemma \ref{lem:some_epoly}, $[\mathcal{O}_{J_\pm}]=\lL^2-1$ and the total contribution is $[Y_{U,U}]=(\lL^2-1)+(\lL^2-1)=2(\lL^2-1)$. 

For even $n$, we have two solutions over each point of $\mathcal{O}_{J_+}$ and no solutions over $\mathcal{O}_{J_-}$, thus the total contribution is again $2(\lL^2-1)$. Observe that the result is independent of the parity of $n$.

We now focus on $X_{T,T}$. By Theorem \ref{tma:covering}, there is a natural regular bijection $\widetilde{X}_{T,T}^0/N_{T,T}\to X_{T,T}$, inducing equality of their relative virtual classes over $\comvar$, where $N_{T,T}=N_G(T)/T\cong \Z_2$. Recall that $X^0_{T,T}=\{(h_1,g_2,h_2)\in X \ | \ (\Stab_G(h_1),\Stab_G(h_1^m))=(T,T) \}$. Thus, $[Y^0_{T,T}]:=\mathcal{Z}(\mathbf{T}^\dagger)([\varphi^0_{T,T}:X^0_{T,T} \to \comvar])$ corresponds to the virtual class of matrices satisfying $h^m=\bar{h}^n=c$ where $c$ is a diagonal matrix with distinct eigenvalues $\lambda,\lambda^{-1}$. This variety can be parametrized by $t\in\C^*$, with $h$ and $\bar h$ having eigenvalues $t^{\pm n}$ and $t^{\pm m}$, respectively, and $c$ having eigenvalues $t^{\pm mn}$. Hence, $t^{2mn}\neq 1$, so the parameter space is $\C^*$ minus the $2mn$-th roots of unity. Thus, in the Grothendieck ring, we have $[Y^0_{T,T}]=\lL-1-2mn$.

We now consider the equivariant virtual classes $[Y^0_{T,T}]^+=[Y^0_{T,T}/\Z_2]$ and $[Y^0_{T,T}]^-=[Y^0_{T,T}]-[Y^0_{T,T}/\Z_2]$, given by the Weyl involution. Under the above parametrization, this involution acts by $t\mapsto t^{-1}$. Its only fixed points in $\C^*$ are $t=1$ and $t=-1$, both of which are excluded from $Y^0_{T,T}$. The remaining $2mn-2$ excluded roots of unity are permuted in pairs, giving $mn-1$ additional points in the quotient. Since $\C^*/\Z_2\cong\mathbb A^1$, the quotient $Y^0_{T,T}/\Z_2$ is isomorphic to $\mathbb A^1$ minus $mn+1$ points. Therefore,
$$
[Y^0_{T,T}]^+=\lL-mn-1, \qquad [Y^0_{T,T}]^-= (\lL-1-2mn)-(\lL-mn-1)=-mn.
$$

We now apply the formula given in Lemma \ref{lem:plusminusformula} (its multiplicativity hypothesis holds since, as in \cite[Example~3.6.22]{vogelthesis}, $Y^0_{T,T}$ is a projective line minus finitely many points, so its class belongs to $V_{\mathbb Z_2}$), which gives $[Y_{T,T}]=[\SL_2/T \times Y^0_{T,T}]^+=[\SL_2/T]^+[Y^0_{T,T}]^+ + [\SL_2/T]^-[Y^0_{T,T}]^-$. Substituting the values of $[\SL_2/T]^\pm$ from Lemma \ref{lem:some_epoly} and those of $[Y^0_{T,T}]^\pm$ computed above, we obtain
\[
[Y_{T,T}]=\lL^3-(mn+1)\lL^2-mn\lL.
\]
To conclude, we determine  $[\widetilde{Y}^0_{T,T}]:=\mathcal{Z}(\mathbf{T}^\dagger)([\widetilde{\varphi}_{T,T}^0:\widetilde{X}_{T,T}^0\to \comvar])$, where, as above, $\widetilde{X}_{T,T}^0:=(X^0_{T,T}\times \SL_2) / T$. We remark that $T$ acts trivially on $X_{T,T}^0$, since all elements in this stratum are diagonal. Therefore, the computation of the virtual class is reduced to the product $[\widetilde{Y}_{T,T}^0]=[Y_{T,T}^0]\cdot[ \SL_2 / T]$, so $[\widetilde{Y}_{T,T}^0]=(\lL-1-2mn)(\lL^2+\lL)$, which establishes the result.
\end{proof}

In analogy with the $\AGL_1$ case, we can compute the virtual class of the representation variety with the matrix product $\mathbf{w}\cdot\mathcal{Z}(\Pi_{n,m})^{d-1}\cdot\mathbf{v}$. Indeed, as observed before Proposition \ref{prop:cerrar_sl2}, the stratification $X=X_G\sqcup X_{T,T}\sqcup X_{U,U}$ gives
$$
[\varphi]=[\varphi_G]+[\varphi_{T,T}]+[\varphi_{U,U}],
$$
so that $\mathbf{v}$ is the column vector with entries $1,1,0,1$ in the basis $\{[\varphi_G], [\varphi_{T,T}], [\widetilde{\varphi}^0_{T,T}],[\varphi_{U,U}]\}$. Repeating the calculation as in Theorem \ref{tma:virtual_class_AGL1}, we obtain the following result:

\begin{tma}
    For $G=\SL_2(\C)$ and $d,n,m\geq 1$, with $n$ and $m$ coprime and $m$ odd, the virtual class of the torus link representation variety is
\begin{align*}
[R_{\SL_2(\C)}(L_{n,m}^d)] &= \left(1 + \frac{m-1}{2}(\lL^2+\lL)\right) \left(2 + (n-1)(\lL^2+\lL)\right) (\lL^3-\lL)^{d-1} \\
&\quad + \left(\lL^3 - (mn+1)\lL^2 - mn\lL\right) (\lL+1)^{d-1} \\
&\quad + \frac{1}{2} (\lL - 1 - 2mn)(\lL^2+\lL) \left( (\lL-1)^{d-1} - (\lL+1)^{d-1} \right) \\
&\quad + 2(\lL^2-1)(2\lL)^{d-1}.
\end{align*}
\end{tma}

\begin{obs}
Since $n$ and $m$ are coprime, at least one of them is odd. Moreover, the right-hand side of the previous formula is symmetric in $n$ and $m$. Therefore, as $L^d_{n,m}$ and $L^d_{m,n}$ are isotopic, the formula holds without the parity assumption on $m$. Note also that it has integer coefficients, since $(\lL-1)^{d-1} - (\lL+1)^{d-1}$ has even coefficients and $(m-1)(n-1)$ is even.
\end{obs}

As a direct application of the TQFT, consider the Hopf link with $n$ twists, denoted by $H_n=L_{n,1}^2$. Applying the preceding theorem with $d=2$ and $(n,m)=(n,1)$, we recover the formula for the virtual class of its $\SL_2(\C)$-representation variety, whose $E$-polynomial was obtained in \cite{gonzalez2023hopflinks}, as stated in the following corollary.

\begin{coro}
    The virtual class of the $\SL_2(\C)$-representation variety of the Hopf link with $n$ twists, $H_n$, is given by
    \begin{equation*}
        [R_{\SL_2(\C)}(H_n)]=((n-1)\lL^2+n\lL-n+5)(\lL^3-\lL).
    \end{equation*}
\end{coro}

\subsection*{AI Tool Use} 
The authors used ChatGPT 6 Astra (OpenAI) and Claude Opus 5.5 (Anthropic) for manuscript editing and proofreading, including LaTeX typesetting, TikZ diagrams, grammatical improvements, and bibliography formatting. All aspects of the research design, formulation of the mathematical results, and their proofs were carried out by the authors, who take full responsibility for the content of the work. No funding from AI companies was received in connection with this research.

\bibliographystyle{abbrv}
\bibliography{mybibliographytoruslinks}
\end{document}